\documentclass{article}       

\usepackage{amsmath,amsfonts}

\def\be{\begin{equation}}
\def\ee{\end{equation}}

\usepackage{graphicx}
\begin{document}

\title{Interpolatory quadrature rules for oscillatory integrals\thanks{Veerle Ledoux is a Postdoctoral Fellow of the Research Foundation-Flanders (FWO).}
}


\author{Veerle Ledoux         \and
        Marnix Van Daele 
}



\maketitle

\newtheorem{myAlgo}{Algorithm}

\begin{abstract}
In this paper we revisit some quadrature methods for highly oscillatory integrals of the form $\int_{-1}^1f(x)e^{{\rm i}\omega x}dx, \omega>0$. Exponentially Fitted (EF) rules depend on  frequency dependent nodes which start off at the Gauss-Legendre nodes when the frequency is zero and end up at the endpoints of the integral
when the frequency tends to infinity. This makes the rules well suited for small as well as for large frequencies. However, the computation of the EF nodes is expensive due to iteration and ill-conditioning. This issue can be resolved by making the connection with Filon-type rules. By introducing some $S$-shaped functions, we show how Gauss-type rules with frequency dependent nodes can be constructed, which have an optimal asymptotic rate of decay of the error with increasing frequency and which are effective also for small or moderate frequencies. These frequency-dependent nodes can also be included into Filon-Clenshaw-Curtis rules to form a class of methods which is particularly well suited to be implemented in an automatic software package.
\end{abstract}

\section{Introduction}
The problem of evaluating oscillatory integrals arises in many applications and forms one of the main classes of integrals for which conventional methods, such as Gaussian or Clenshaw-Curtis quadratures, are inadequate. The integrands of interest to us are the ones which typically appear in  the finite Fourier transform
\begin{equation}I[f]=\int_{-1}^1 f(x)e^{{\rm i}\omega x}dx.\label{prob}\end{equation}
The evaluation of \eqref{prob} is a problem that occurs frequently throughout the applied sciences. The availability of efficient highly oscillatory quadrature rules for \eqref{prob} is e.g.\ essential to a number of recent methods for highly oscillatory differential equations (see \cite{Degani2006,Iserles2002,Iserles2004N,Ledoux2010}) or in the implementation of boundary integral equation methods in high-frequency scattering \cite{Dominguez}. 
Note that studying \eqref{prob} covers the more general integration interval $[a,b]$ by the following identity
\begin{equation}
\int_a^b f(x)e^{i\omega x}dx=\frac{b-a}{2}e^{i\omega(b+a)/2}\int_{-1}^1 g(t)e^{i {\hat \omega} t }dt
\end{equation}
with ${\hat \omega}=(b-a)\omega/2$ and $g(t)=f((b-a)t/2+(b+a)/2)$.

\begin{figure}
	\centering
	\includegraphics[width=0.85\textwidth]{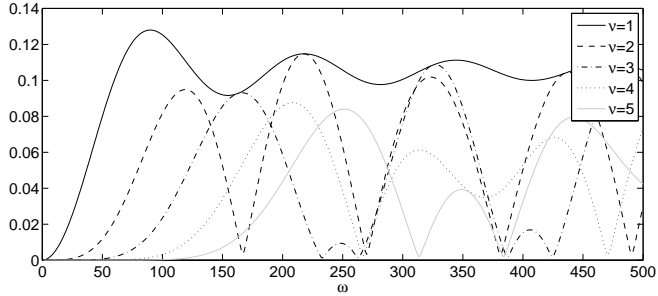}
	\caption{The error in Gauss-Legendre quadrature in approximating $\int_0^{1/10} e^x e^{{\rm i}\omega x}dx$ for different values of $\omega$. }
		\label{fig:fig0}
\end{figure}

When the frequency parameter $\omega$ in \eqref{prob} is large, the integrand is highly oscillatory. In this case, a prohibitively large number of quadrature nodes is needed if one uses a quadrature rule based on polynomial interpolation of the integrand such as the standard Gaussian quadrature rule
\begin{equation}\label{gauss}
I[f]\approx Q^G_\nu[f]=  \sum_{l=1}^\nu b_l f(c_l) e^{{\rm i}\omega c_l}\end{equation}
where $c_1,\dots,c_\nu \in [-1,1]$ are distinct nodes, and $b_1,\dots,b_\nu$ are interpolatory weights. If the nodes are selected as the zeros of the Legendre polynomial of degree $\nu$, then the quadrature rule \eqref{gauss} is of the order $2\nu$. By its very construction, the Gauss-Legendre rule is exact if the integrand is a polynomial of degree ${2\nu-1}$. However, if the integrand oscillates rapidly, and unless we use a huge number of function evaluations, the polynomial interpolation underlying the classical Gauss rule is useless. This is illustrated in Figure \ref{fig:fig0}. We have computed
\begin{equation}\int_0^{1/10} e^x e^{{\rm i}\omega x}dx=\frac{-1+e^{(1+{\rm i}\omega)/10}}{1+{\rm i}\omega},\label{prob10}\end{equation}
by Gaussian quadrature with different numbers of nodes. The figure displays the absolute value of the error as a function of the frequency $\omega$. For small $\omega$ the standard Gaussian quadrature rules do fine, but for larger $\omega$ high oscillation sets in soon and the error hardly reduces when the number of quadrature points increases.

The accuracy of the approximation \eqref{gauss} can be greatly improved by making the weights $b_l$ and quadrature nodes $c_l$ $\omega$-dependent. In fact, methods can be constructed which have a performance that drastically improves as the frequency grows. A first class of such methods tuned to oscillatory integrands is formed by applying a technique called exponential fitting (EF) to the Gauss rules (see \cite{Ixarubook}). The Filon-type rules form a second class of quadrature rules with frequency dependent nodes. These rules are based on replacing $f(x)$ with an interpolating polynomial. The paper \cite{Iserles2004} analysed some particular types of the Filon approach and in \cite{Iserles2005} it was shown that the asymptotic rate of decay of the error with increasing frequency can be improved by also interpolating the derivatives of $f$ at the endpoints of the interval. Alternatively, derivatives can be avoided by allowing the interpolation nodes to approach the endpoints as $\omega$ increases, see \cite{Iserles2004b}. 

In sections 2 and 3 we give the basic ideas of the exponential fitted Gauss methods and the Filon methods. In section 4 the connection between the two classes of methods is illustrated by considering the construction of some numerical schemes with different numbers of quadrature nodes. We show how to form quadrature rules which combine the strengths from both the EF and the Filon technique. They have optimal asymptotic behaviour and reduce to the classical Gauss methods for small frequencies $\omega$, while avoiding the iterative process which is needed to compute the EF nodes. In literature \cite{Iserles2004b,Iserles2005,Deano2009,Huybrechs2009}, much attention has been devoted to the improvement of the asymptotic convergence rates of Filon-type methods but little has been written on how to exploit the high asymtotic orders in a practical implementation of the schemes. In fact, these methods are often not very well suited to be implemented in a ready-to-use computer code where an adaptive procedure is used to reach a user-requested accuracy. Moreover, numerical difficulties appear, putting a `practical' limit on the highest asymptotic order achievable.
In section 5, we examine the addition of Chebyshev nodes to the methods of section 4. This allows us to form practical schemes which are well suited to be applied in an adaptive procedure and thus in an automatic software package.

\section{Exponentially fitted Gauss rules}
A first approach towards the determination of $\omega$ dependent parameters for the Gauss quadrature rule consists in the use of the exponential fitting techniques of Ixaru et al. \cite{Ixaru2001,Ixarubook,Kim2003,VanDaele2005}.
Where a classical approximation formula is designed to be exact for polynomials of low degree, an exponentially fitted method is designed to be exact when the integrand is some suitably chosen combination of exponential functions, perhaps with polynomial terms, or products of polynomials and exponentials.

In \cite{VanDaele2005}, a variety of exponentially fitted Gaussian rules was discussed.
The authors discuss the construction of a $\nu$-point quadrature formula $\sum_{k=1}^\nu w_k l(x_k)$ for 
\begin{equation}\int_{-1}^1 l(x)dx,\label{eq0ef}\end{equation}
where $l(x)$ shows an oscillatory behaviour with frequency $\omega$. 
The functional 
\[{\mathcal L}[l;x;h;{\bf a}]=\int_{x-h}^{x+h}l(z)dz-h\sum_{k=1}^\nu w_k l(x+{c}_kh),\quad {c}_k\in[-1,1]\]
is considered, where ${\bf a}=[{c}_1,{c}_2,\dots,{c}_\nu, w_1,w_2,\dots,w_\nu]$ is a vector of $2\nu$ unknowns which can all depend on ${\omega}$. To make the connection with the problem \eqref{eq0ef} $x=0$, $h=1$ will be taken. 
The following (reference) set of functions $\ell(x)$ is considered:
\begin{align}\label{set}
&1,x,x^2,...x^K, \\
\label{set1}
&\exp(\pm{\bar \omega}x), x\exp(\pm{\bar \omega}x), x^2\exp(\pm{\bar \omega}x),\dots,x^P\exp(\pm{\bar \omega x}),
\end{align}
where $K+2(P+1)+1=2\nu$. 
The reference set for a $\nu$-point rule is thus characterized by two integers: $K$ and $P$. The parameter $P$ is called the level of tuning. The set in which there is no exponential fitting component, corresponding to a traditional Gauss method is identified by $P=-1$, while the set in which there is no purely polynomial component is identified by $P=\nu-1$. As explained in \cite{Ixaru2001}, one can either consider the set of power functions \eqref{set} or the exponential fitting set \eqref{set1} and compute ${\mathcal L}$ for each of these functions. The parameters in ${\bf a}$ can then be fixed by imposing the condition that as many terms as possible in the expressions for ${\mathcal L}$ vanish. The set of functions as shown in \eqref{set1} covers both the oscillatory and the exponential case, but since we consider only oscillatory integrals here, we take ${\bar \omega}=i\omega$. 

In \cite{VanDaele2005} it was shown that the pure exponentially fitted case ($P=\nu-1$) forms the best option for highly oscillatory integrals.
When symmetric weights and antisymmetric abscissae are assumed, the following system needs to be solved to construct a pure EF rule:
\be H_m(\omega,{\bf {c}},{\bf w})=0,\quad m=0,...,P,\label{sysH}\ee
with
\begin{eqnarray}
H_0(\omega,{\bf {c}},{\bf w})&=&2\sin\omega-\omega\sum_{k=1}^\nu w_k \cos(\omega {c}_k),\\
H_i(\omega,{\bf {c}},{\bf w})&=&\frac{\partial}{\partial \omega}H_{i-1}(\omega,{\bf {c}},{\bf w})=\frac{\partial^i}{\partial \omega^i}H_{0}(\omega,{\bf {c}},{\bf w}),\quad i=1,2,\dots.
\end{eqnarray}
To find the solution to the system \eqref{sysH} an iteration procedure is used. As reported in \cite{Ixaru2001} some ill conditioning issues appear for large $\omega$.
\begin{figure}
	\centering
		\includegraphics[width=0.75\textwidth]{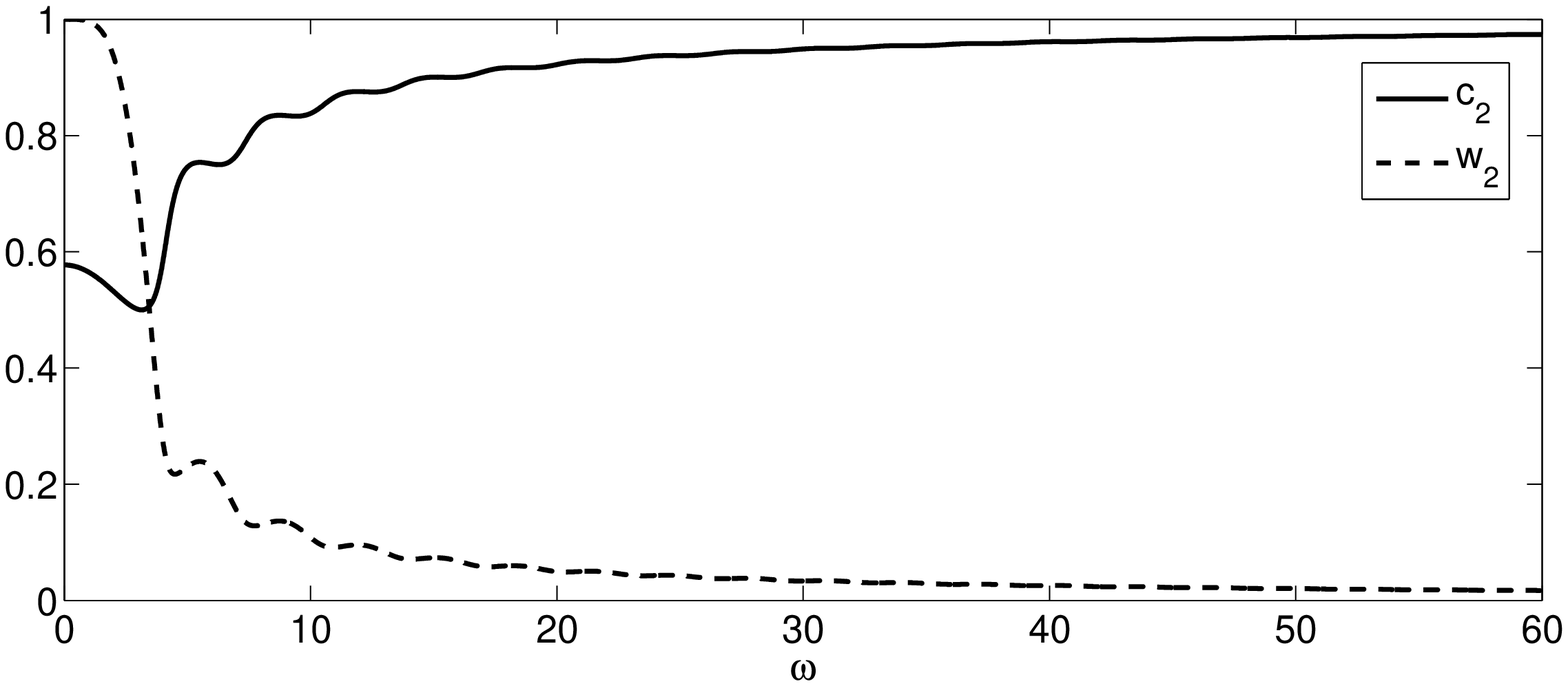}
		\includegraphics[width=0.75\textwidth]{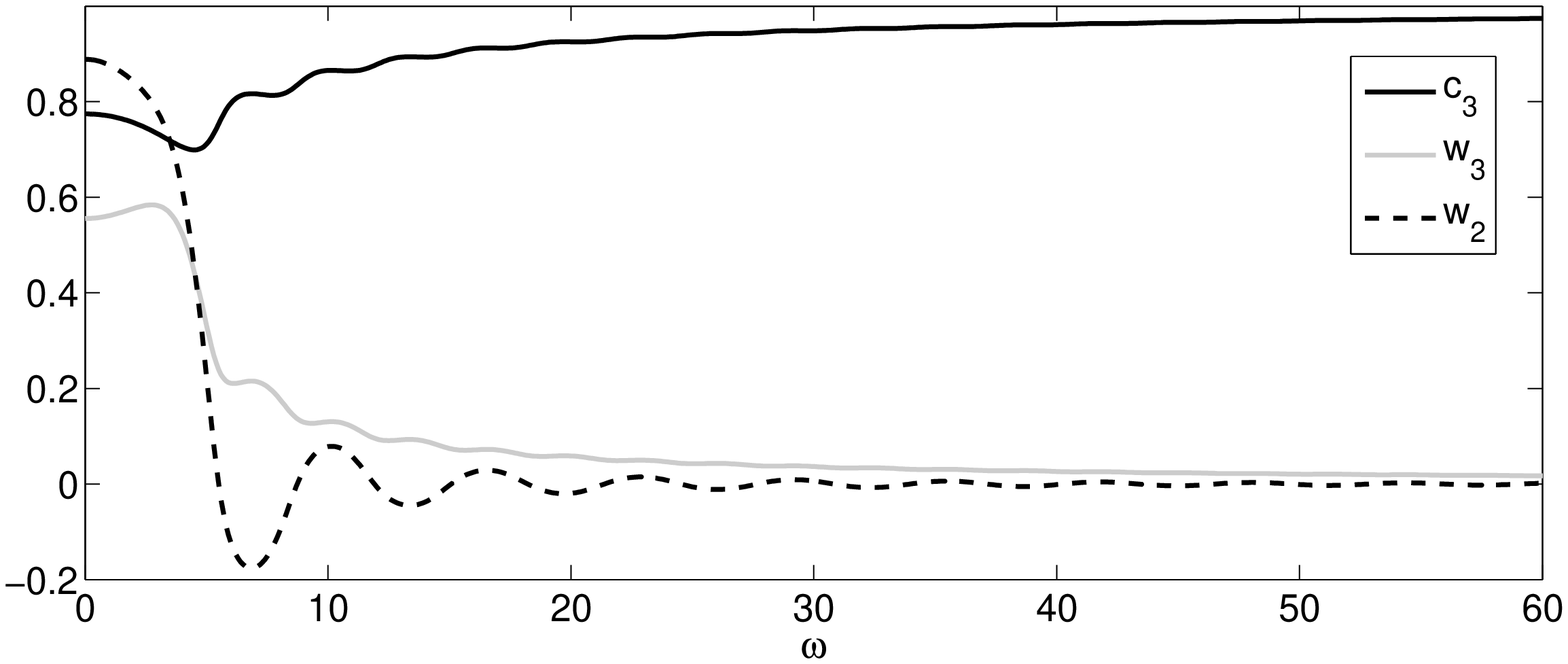}
		\includegraphics[width=0.75\textwidth]{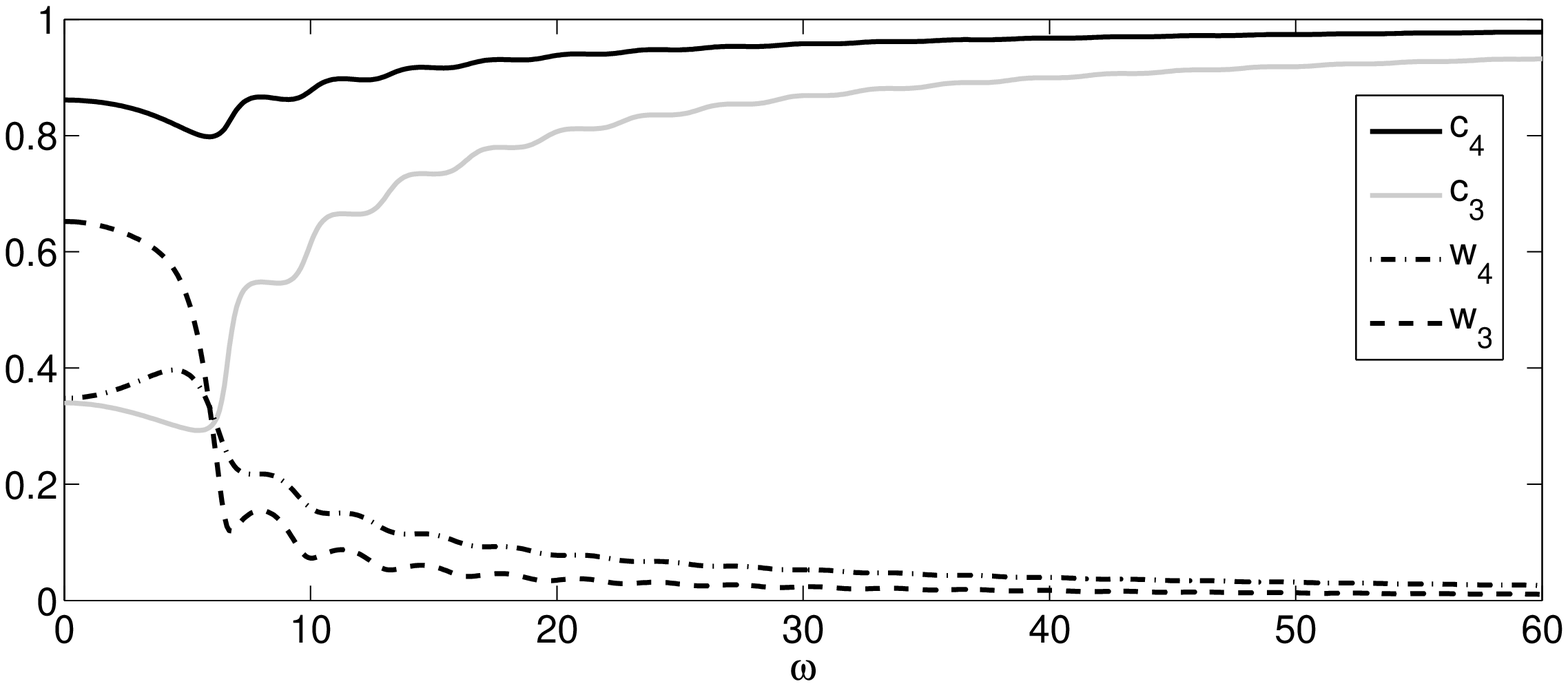}
	\caption{Nodes and weights of some EF rules with $\nu$ nodes: (1) the ${c}_2(\omega)$ and $w_2(\omega)$ curve for $\nu=2$, (2) The ${c}_3(\omega)=-{c}_1(\omega)$, ${w}_1(\omega)={w}_3(\omega)$ and $w_2(\omega)$ curves for $\nu=3$, and (3) ${c}_3(\omega),{c}_4(\omega) $ and ${w}_4(\omega),{w}_3(\omega)$ for $\nu=4$ .}	\label{fig:figef2b}
\end{figure}
Figure \ref{fig:figef2b} shows the nodes and weights for some EF-schemes.
 All EF rules reduce to the classical $\nu$-point Gauss-Legendre method in the limiting case $\omega=0$ and consequently have a similar error behaviour for small $\omega$. 
In \cite{Coleman2006}, it was shown that for larger $\omega$, the quadrature error decays like $\omega^{{\bar \nu}-{\nu}}$, with ${\bar \nu}=\left\lfloor (\nu-1)/2\right\rfloor$ as $\omega \to \infty$. In particular for $\nu=1$ the error decays as $\omega^{-1}$, for $\nu=2,3$ we have an $\omega^{-2}$ error, but for $\nu=4,5$ we have an even faster decay rate of $\omega^{-3}$.

\section{Filon-type methods}
A method which recently gained popularity in the context of highly oscillatory quadrature is the Filon method, based on an idea originally due to Filon \cite{Filon}.
Where classical Gauss quadrature interpolates the whole integrand $f(x)e^{{\rm i}\omega x}$ at distinct nodes $c_1<c_2<\dots<c_{\nu}$ in $[-1,1]$ by a polynomial $p$ of degree ${\nu -1}$ and defines the approximation as the exact evaluation of the result,
we now interpolate only the function $f(x)$ at $c_1,\dots,c_\nu $ by a polynomial ${\bar f}(x)$ and define the Filon approximation as the exact evaluation of this result, i.e.
\[I[f]\approx Q_\nu^F[f]=\int_{-1}^1 {\bar f}(x)e^{{\rm i}\omega x}dx= \sum_{l=1}^\nu b_l(\omega)f(c_l)\]
where $b_l(\omega)=\int_{-1}^1\ell_l(x) e^{{\rm i}\omega x}dx$ and  $\ell_l$ is the $l$th cardinal polynomial of Lagrangian interpolation. 

For small $\omega$, the Filon-type method has the same `classical' order as the corresponding traditional Gauss method with the same quadrature nodes $c_1<c_2<\dots<c_{\nu}$ (see \cite{Iserles2004}).
For an increased frequency $\omega$, a Filon-type method typically results in smaller errors and this behaviour is enhanced once the endpoints are included among the quadrature nodes.
This behaviour can be explained using the asymptotic expansion of the integral considered. As shown in \cite{Iserles2004b,Iserles2005}, an explicit asymptotic expansion of 
\[I[f]=\int_{-1}^1 f(x)e^{{\rm i}\omega x}dx\] 
can be derived by repeated integration by parts:
\begin{equation}\label{asy}
I[f]\approx -\sum_{m=0}^\infty \frac{1}{(-{\rm i}\omega)^{m+1}}\left[e^{{\rm i}\omega}f^{(m)}(1)-e^{-{\rm i}\omega }f^{(m)}(-1)\right].
\end{equation}
If we define (as in \cite{Iserles2004b}) $v={\bar f}-f$ and use \eqref{asy}, we have
\begin{align}
\nonumber Q^{F}_{\nu}[f]-I[f]&=I[{\bar f}]-I[f]=I[v]\\
&\sim  -\sum_{m=0}^\infty \frac{1}{(-{\rm i}\omega)^{m+1}}\left[e^{{\rm i}\omega }v^{(m)}(1)-e^{-{\rm i}\omega }v^{(m)}(-1)\right],\quad |\omega|\gg 1.\label{qefnu}
\end{align}
When $c_1=-1$ and $c_\nu=1$, it follows from interpolation conditions that $v(-1)=v(1)=0$, which gives $Q_\nu^F[f]-I[f]=O(\omega^{-2})$.
It is clear, however, from \eqref{asy} that the asymptotic estimate can be further improved by letting the interpolating polynomial ${\bar f}$ depend on derivatives of $f$. This method was proposed in \cite{Iserles2005}. Instead of Lagrange interpolation, Hermite interpolation is then used. In this case $Q_\nu^F[f]-I[f]=O(\omega^{-p-1})$, i.e. an `asymptotic order' $p+1$ can be reached where $p$ is the number of derivatives at the endpoints: ${\bar f}^{(l)}(-1)={f}^{(l)}(-1),{\bar f}^{(l)}(1)={f}^{(l)}(1), l=0,\dots, p-1$.  The same arbitrarily high asymptotic order can also be achieved without computation of derivatives (which is often expensive) but by allowing the interpolation points to depend on $\omega$, as noted in \cite{Iserles2004b}. Instead of derivatives, finite difference approximations are then used, with spacing inversely proportional to the frequency leading to a new family of Filon-type methods, called the {\em adaptive} Filon-type method. 
We will see in further sections that there is a close connection between (pure) EF methods and adaptive Filon methods, both using $\omega$ dependent interpolation nodes to reach an optimal asymptotic error for a particular number of nodes $\nu$. A pure EF scheme with $\nu$ nodes is exact for $f(x)=1,x,x^2,\dots, x^{\nu-1}$. This means that we can see the application of an EF scheme to a problem of the form \eqref{prob} as a specific Filon-type method. That is, EF involves the replacement of $f$ by ${\bar f}$, a specific polynomial of degree $\nu-1$ through the EF nodes and the computation of the analytic solution of this approximating integral.
 We will exploit this connection between the two approaches to construct methods which share the property of optimal behaviour for both small and large $\omega$ values with the EF rules, while avoiding the need for iteration or the ill-conditioning issues when computing the frequency dependent nodes.

Note that, like the EF rules, the computation of the Filon approximation is based on the ability to compute the moments
\[\int_{-1}^1 x^k e^{{\rm i}\omega x} dx.\]
For this particular oscillator, the moments can be computed in closed form, either through integration by parts or by the identity
\be\int_{-1}^1 x^k e^{{\rm i}\omega x} dx = \frac{1}{(-{\rm i}\omega)^{k+1}}[\Gamma(1+k,{\rm i}\omega)-\Gamma(1+k,-{\rm i}\omega)],\label{Gamma}\ee
where $\Gamma$ is the incomplete Gamma function.

\section{Numerical schemes}\label{sec4}
\subsection{Methods with $\nu=1$ quadrature nodes}
We use the approaches discussed in the previous sections to construct one-node quadrature rules for the problem
\begin{equation}\label{prob2} I[f]=\int_{-1}^1 F(x)dx=\int_{-1}^1 f(x) \exp({\rm i} \omega x)dx.\end{equation}
We will be able to conclude that only one $\nu=1$ scheme is suitable for the integration of \eqref{prob2} when the integrand is highly-oscillatory.
For reasons of symmetry we choose the quadrature node in the middle of the integration interval ($c_1=0$) and we will see that this forms indeed the best choice.

 The EF quadrature rule with one quadrature node ${c}_1=0$, is exact for the integrand $\exp({\rm i}\omega x)$, which means that the relation $\int_{-1}^1 \exp({\rm i}\omega x)dx - w_1 \exp(0)=0$
should be satisfied.
This gives us the weight $w_1=2 \sin(\omega)/\omega$ and the second order EF rule for our problem is 
\be\label{ef1}
Q^{EF}_1[F]=
\frac{2\sin({ \omega})}{{ \omega}} f(0).
\ee

Applying a one-node Filon-type rule consists in replacing $f(x)$ in $I[f]$ by a constant ${\bar f}$ and using the exact integral $I[{\bar f}]$ as approximation for $I[f]$.  If ${\bar f}=f(0)$, the Filon-type method is exactly the same as the EF scheme \eqref{ef1}, i.e.\ $Q^{F}_1[f]=Q^{EF}_1[f]$.

\begin{figure}
	\centering
	\includegraphics[width=0.8\textwidth]{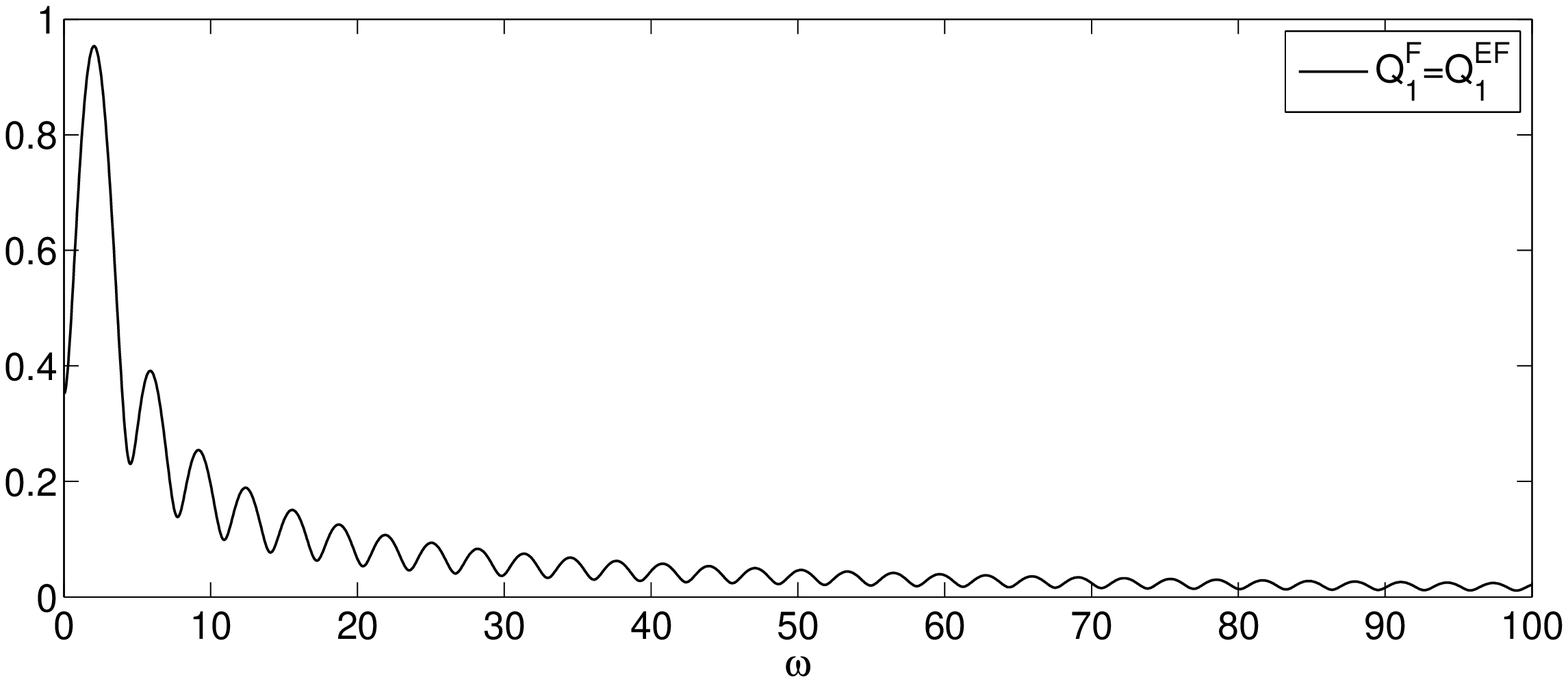}
		\includegraphics[width=0.8\textwidth]{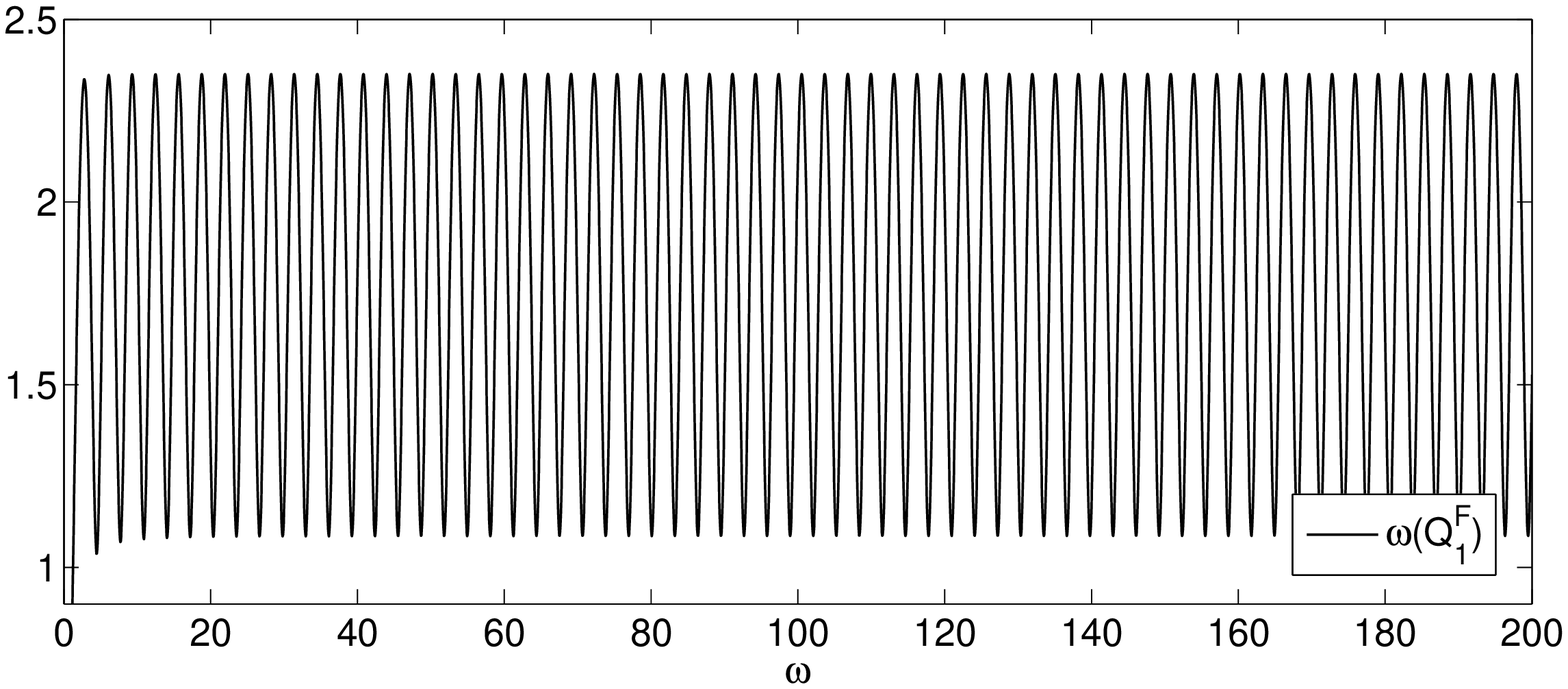}
	\caption{The error in the exponentially fitted (and Filon)  method with $\nu=1$ for  $\int_{-1}^1 e^x e^{{\rm i}\omega x}dx$ and different values of $\omega$. The top graph shows the absolute error $E^{[EF]}(=E^{[F]})$, the bottom the normalised error $\omega E^{[EF]}$.}	\label{fig:fig1}
\end{figure}

In Figure \ref{fig:fig1} we display the error for the problem 
\[\int_{-1}^1 e^x e^{{\rm i}\omega x}dx=\frac{e^{1+{\rm i} \omega}-e^{-(1+{\rm i}\omega)}}{1+{\rm i}\omega}\]
 for this Filon/EF method. The bottom figure confirms that the error decays like $O(\omega^{-1})$. This is the most optimal asymptotic behaviour which can be reached for a method with one quadrature point. It can be shown that for all possible choices (real or complex) of quadrature nodes $c_1$, the same asymptotic order is obtained as for $c_1=0$. For $\omega=k\pi, (k=1,2,\dots)$ for instance, we have $I[{\bar f}]=0$ and an ${\bar f}$-independent quadrature error $I[f]-I[{\bar f}]=I[f]=O(\omega^{-1})$. Since the quadrature node $c_1=0$ also optimalizes the error behaviour for smaller $\omega$ values, $f(0)$ is a good choice for ${\bar f}$.

\subsection{Methods with $\nu=2$ quadrature nodes}
\begin{figure}
	\centering
		\includegraphics[width=3.9cm]{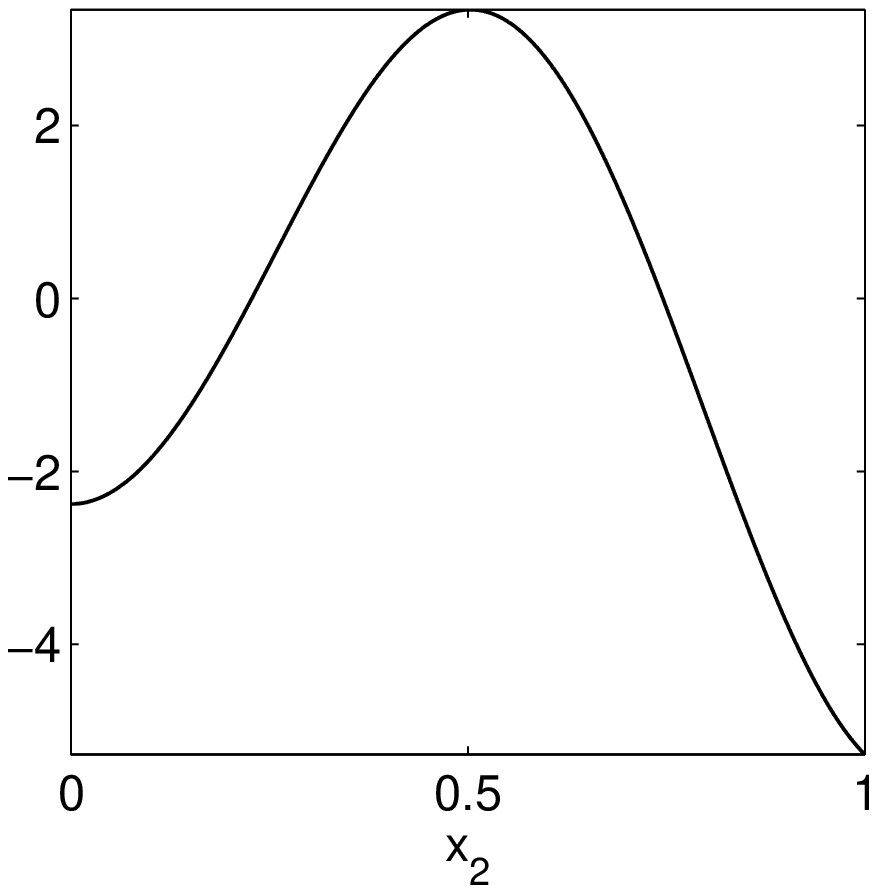}
		\includegraphics[width=3.9cm]{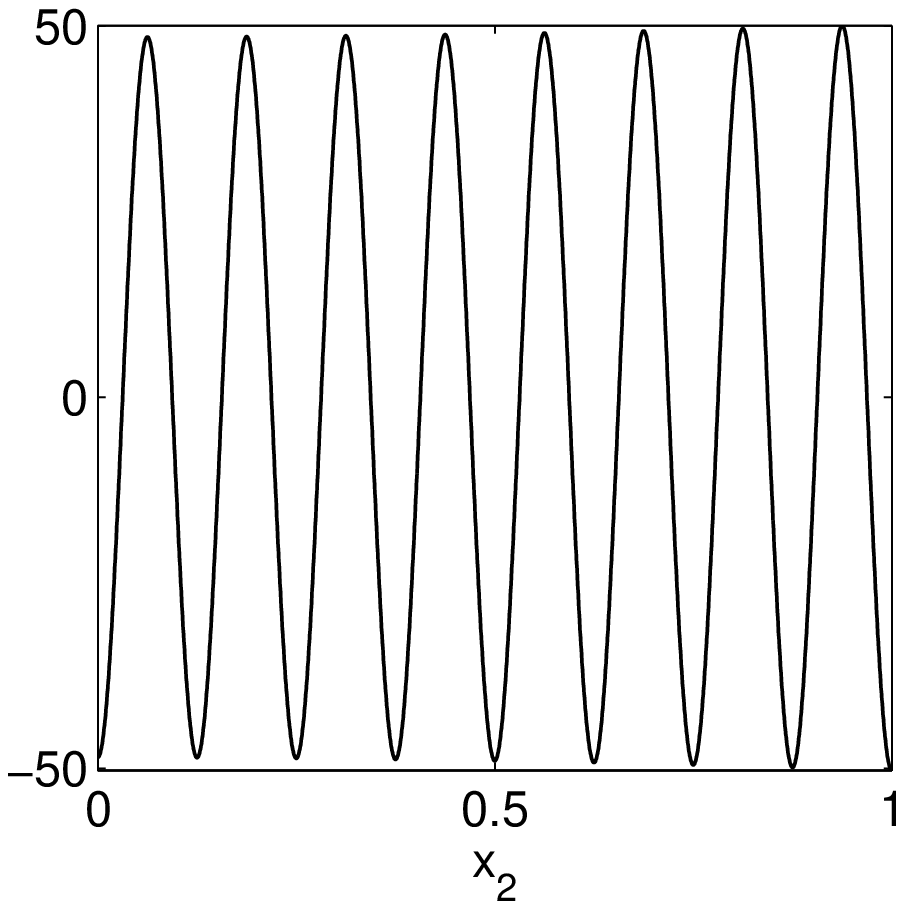}
		\includegraphics[width=3.9cm]{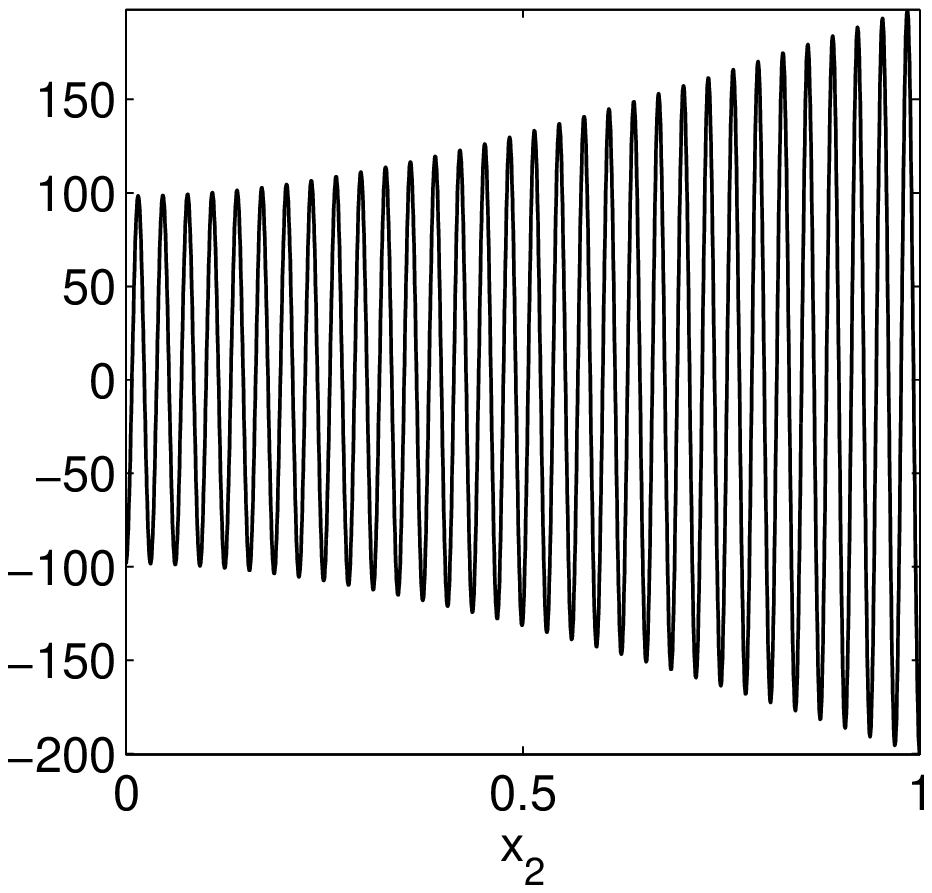}
	\caption{$G(x_2)=(\sin\omega-\omega\cos\omega)\cos(\omega x_2)-\omega x_2\sin\omega \sin(\omega x_2)$ for three different $\omega$-values: $\omega=5$, $\omega=50$ and $\omega=200$.}	\label{fig:fig2x}
\end{figure}

For the EF method, we assume antisymmetric quadrature nodes and symmetric weights, i.e.\ we take ${c}_1=-{c}_2 \in [-1,1]$ and $w_1=w_2$.
The two equations
\[\int_{-1}^1 \exp({\rm i}\omega x)dx-w_2 \exp(-{\rm i}\omega {c}_2)-w_2 \exp({\rm i}\omega {c}_2)=0\]
\[\int_{-1}^1 x \exp({\rm i}\omega x)dx+w_2 {c}_2\exp(-{\rm i}\omega {c}_2)-w_2 {c}_2\exp({\rm i}\omega {c}_2)=0\]
lead to the following relations for $w_2$ and ${c}_2$:
\begin{align}\label{eq0a}
&w_2\omega\cos(\omega {c}_2)-\sin(\omega)=0,\\
&w_2{c}_2\omega^2\sin(\omega {c}_2)-\sin(\omega)+\omega\cos(\omega)=0.\label{eq1a}
\end{align}
As described in \cite{VanDaele2005}, the system \eqref{eq0a}-\eqref{eq1a} allows us to obtain suitable values for ${c}_2$ en $w_2$ for a particular $\omega$.
The EF method then takes the form
\[ Q_2^{EF}[F]=w_2F\left({c}_2\right)+w_2F\left(-{c}_2\right).\]
If $\cos(\omega {c}_2)\neq 0$ then we know from \eqref{eq0a} that $w_2=\sin\omega/(\omega\cos(\omega {c}_2))$. Substituting this into the second relation gives
\begin{equation}
(\sin\omega-\omega\cos\omega)\cos(\omega {c}_2)-\omega {c}_2\sin\omega \sin(\omega {c}_2)=0
\label{xw2}
\end{equation}
which allows us to compute ${c}_2$ for a given $\omega$ value. There is however no unique solution to this equation. As shown in \cite{VanDaele2005}, the best option is to use the so-called classical solution curve which passes through the Legendre point $1/\sqrt{3}$. This curve is shown in Figure \ref{fig:figef2b}.
One observes that the curve for ${c}_2$ starts off at the Legendre node for $\omega=0$ and tends to 1 for $\omega\to\infty$, the weight $w_2$ on the other hand decreases to zero for large $\omega$.
As described in \cite{Ixaru2001}, ${c}_2$ and $w_2$ can be computed from the system \eqref{eq0a}-\eqref{eq1a} using a Newton iteration process. As reported earlier, there are however some ill-conditioning problems. Figure \ref{fig:fig2x} shows $G(x_2)=(\sin\omega-\omega\cos\omega)\cos(\omega x_2)-\omega x_2\sin\omega \sin(\omega x_2)$ (i.e.\ the left-hand side of eq.\ \eqref{xw2}) for three different values of $\omega$. It is clear that $G(x_2)$ vanishes for several $x_2$ values. From the asymptotic expansion \eqref{asy}, we know that for large frequencies the most important piece of information is situated around the endpoints. The ``optimal'' ${c}_2$ value is, consequently, the root of $G(x_2)$ closest to the endpoint 1 (this corresponds to the classical solution curve) . The oscillatory character of $G(x_2)$ increases however with $\omega$ and it is easy to see intuitively that this root is more difficult to locate for high $\omega$. This rootfinding process is avoided when using a Filon approach.

A Filon method replaces $f$  by a first degree polynomial. Consequently the Filon method will be exact for $f(x)=1$ and $f(x)=x$, that is for  $F(x)= \exp ({\rm i}\omega x)$ and $F(x)= x\exp({\rm i}\omega x)$, but also for $F(x)= \exp (-{\rm i}\omega x)$ and $F(x)= x\exp(-{\rm i}\omega x)$. There is thus a close connection with the EF rule. Indeed  the Filon-method and the EF rule coincide when both schemes use the same quadrature nodes, i.e.\ the nodes given by the EF-relations \eqref{eq0a}-\eqref{eq1a}. However, a Filon method can also take other quadrature points, i.e. a different first degree polynomial approximation for $f$.  Standard choices are Legendre nodes or Lobatto nodes. The Legendre nodes lead to a method of classical order 4 and asymptotic order 1 
and the Filon-Lobatto method
 has classical order 2 and asymptotic order 2. The $\omega$ dependent nodes from the EF scheme, starting at the Legendre nodes and tending to the Lobatto nodes for large $\omega$, allow to combine the strengths from both, leading to a method with classical order 4 and asymptotic order 2 (see Figure \ref{fig:filon2}). 
 Since the EF nodes can be difficult and expensive to compute for larger $\omega$, we investigate if other frequency dependent nodes can be used to construct Filon-type methods with the same or even better properties. The nodes are in fact parameters which can be chosen in such a way that the Filon error (both for small and large $\omega$ values) is optimalised. We constructed a (new) {\em adaptive} Filon method with interpolation points  which take the form of the Legendre interpolation points ${\bar c_1}=-1/\sqrt{3}, {\bar c_2}=1/\sqrt{3}$ for $\omega\to 0$ and tend to the endpoints of the integration interval for  $\omega\gg 1$. A good error behaviour for small frequencies and an $\omega^{-2}$ asymptotic behaviour can be combined when the nodes $c_1(\omega)$ and $c_2(\omega)$ remain close to the Legendre nodes ${\bar c_1}$ and ${\bar c_2}$ for smaller frequencies and have an $\omega^{-1}$ attraction towards the endpoints $-1$ and $1$ for large $\omega$ values. The EF nodes satisfy these requirements. 
 
\begin{figure}
	\centering
		\includegraphics[width=0.8\textwidth]{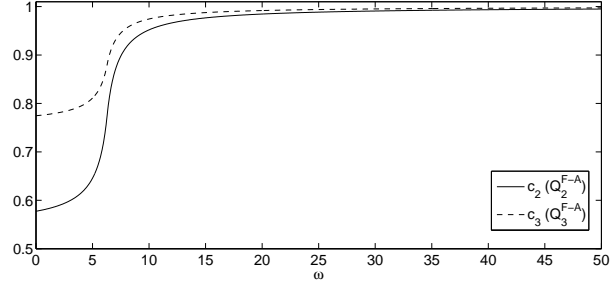}
		\caption{$c_2(\omega)$ of the adaptive Filon method $Q_2^{F\text{-}A}$ with two quadrature nodes and $c_3(\omega)$ of the adaptive Filon method $Q_3^{F\text{-}A}$ with three quadrature nodes}
	\label{fig:figc1}
\end{figure}

\begin{figure}
	\centering
			\includegraphics[width=0.8\textwidth]{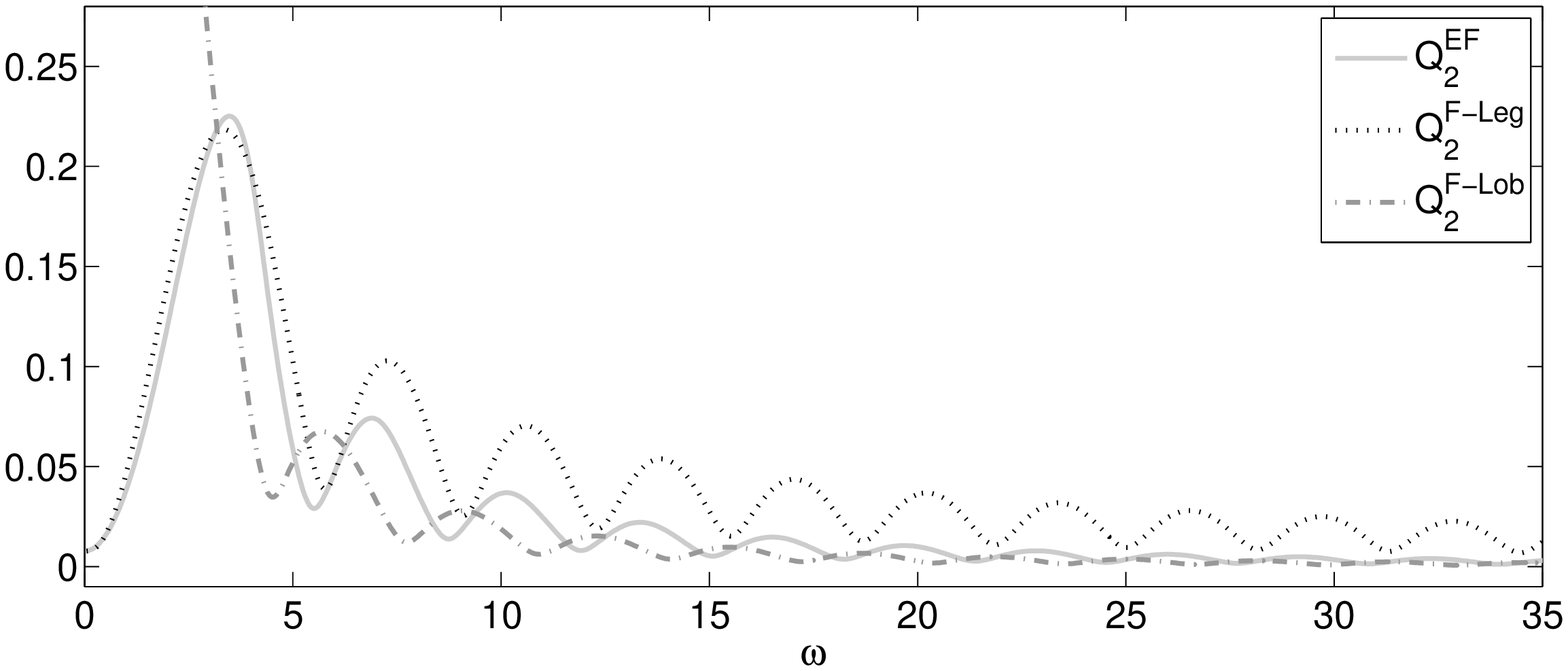}
		\includegraphics[width=0.8\textwidth]{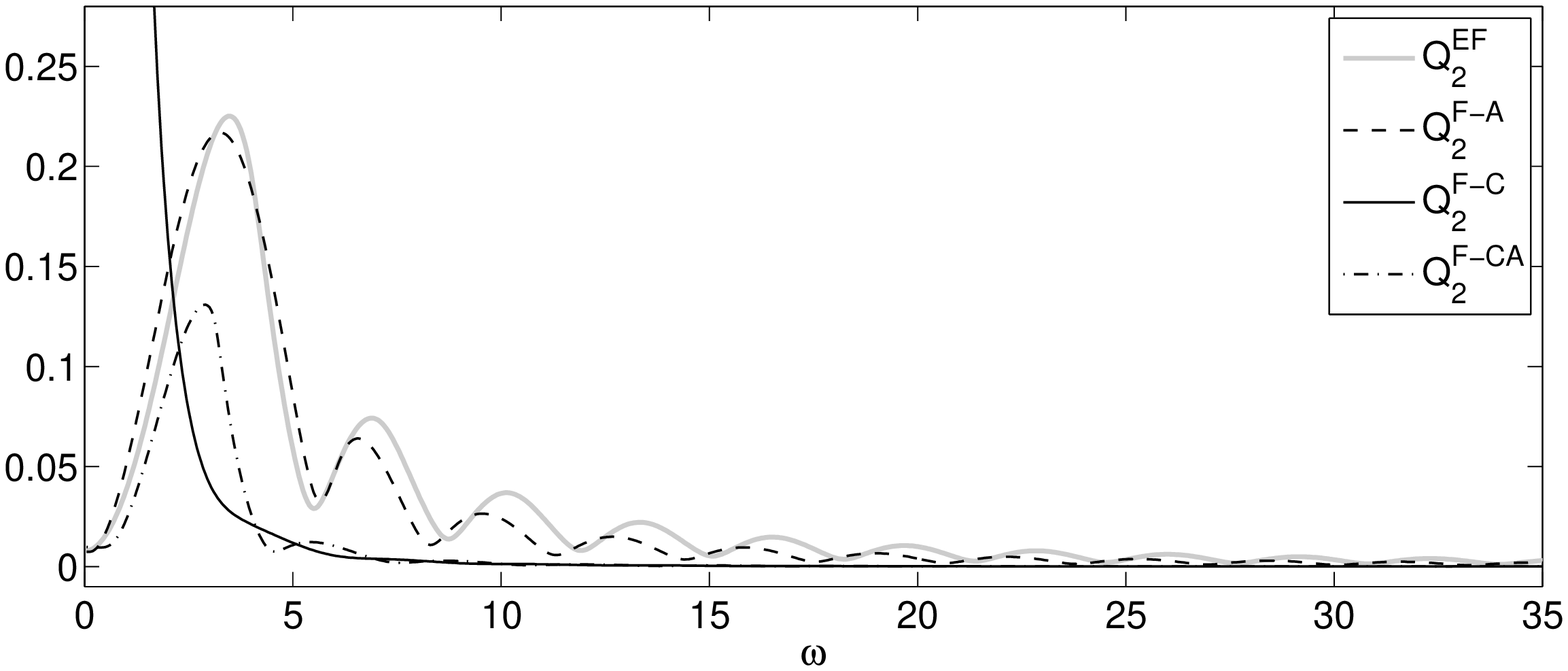}
		\includegraphics[width=0.8\textwidth]{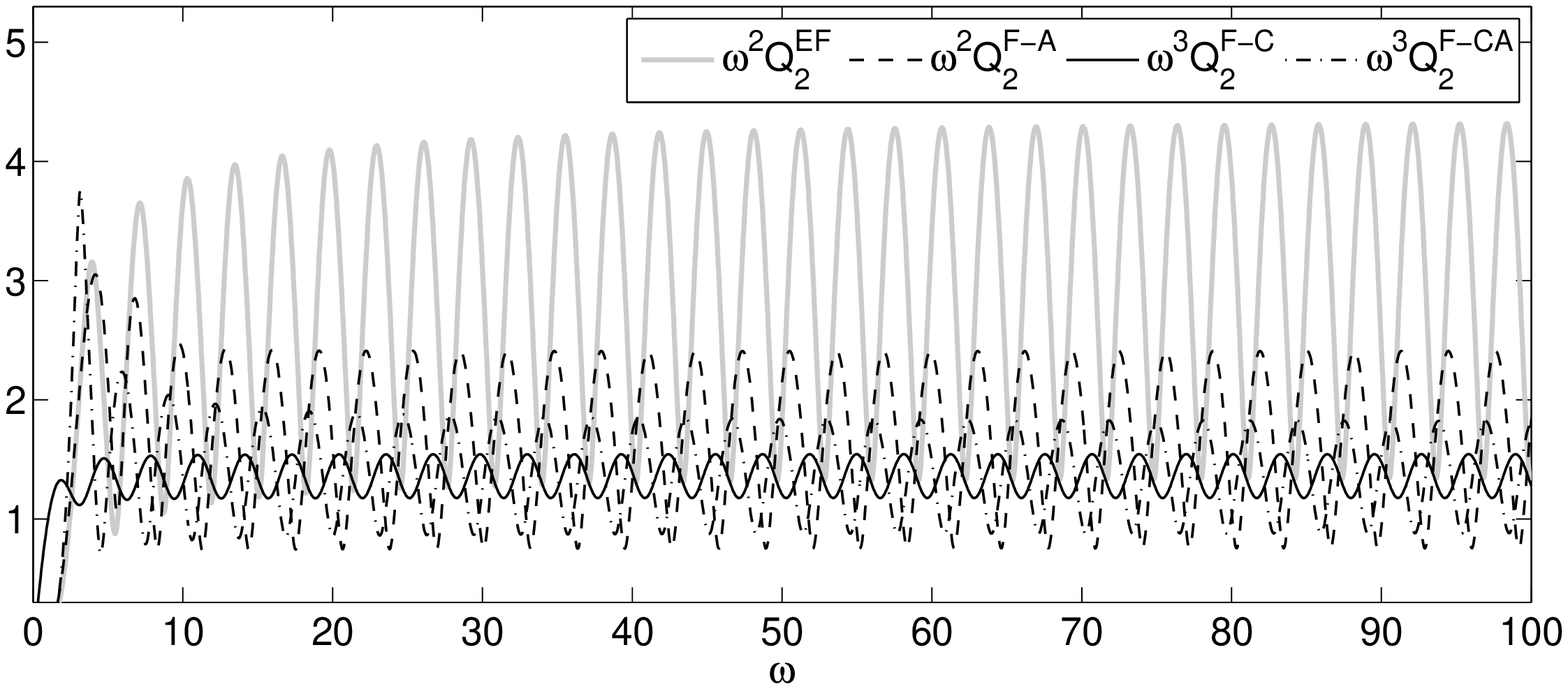}
		\caption{The error in some $\nu=2$ Filon-type schemes for $\int_{-1}^1 e^x e^{{\rm i}\omega x}dx$ and different values of $\omega$. The top graphs show the absolute error and the bottom figure contains the normalised errors.}
	\label{fig:filon2}
\end{figure}

We use the $S$-shaped function $S(\omega;r;1)$ which is centred around $r=2\pi$ to construct a suitable curve for the $c_1$ and $c_2$ nodes: 
\begin{equation}\label{adfilon}
c_2(\omega)=-c_1(\omega)=1-(1-{\bar c}_2)S(\omega;2\pi;1),\quad S(\omega;r;n)=\frac{1-\frac{\omega^n-r^n}{1+|\omega^n-r^n|}}{1+\frac{r^n}{1+r^n}}.
\end{equation}
The interpolation point $c_2(\omega)$ is shown in Figure \ref{fig:figc1}. The point $r=2\pi$ is approximately where the $\omega^{-1}$ decrease to the endpoints starts in the curves for the EF nodes. Results for this adaptive Filon method $Q_2^{F-A}$ are shown in Figure \ref{fig:filon2}. 

The choice \eqref{adfilon} for the nodes of the adaptive Filon method can also be motivated theoretically. Using the asymptotic expansion \eqref{asy} one can show that when one assumes  that the quadrature nodes satisfy the relation $c_2=-c_1$, $c_2$ needs to behave asymptotically as $1+\lambda\omega^{-p},p\geq1$ where $\lambda$ is independent of $\omega$, in order to have asymptotic order 2. 
Let us define $v={\bar f}-f$ and use the expansion \eqref{asy} to obtain an asymptotic expression for the error of the Filon method:
\be 
Q_2^{F}[f]-I[f]=I[v]\sim  -\sum_{m=0}^\infty \frac{1}{(-{\rm i}\omega)^{m+1}}\left[e^{{\rm i}\omega}v^{(m)}(1)-e^{-{\rm i}\omega}v^{(m)}(-1)\right],\quad |\omega|\gg 1.\label{asyFilon}
\ee
We have one interpolation node ${c}_1=-1+\sigma_1$ in the vicinity of the endpoint $-1$ and one node ${c}_2=1+\sigma_2$ in the vicinity of the endpoint $1$. The interpolation error in the vicinity of $1$ is then given by 
\be v(x)=s_1(x) (x-1-\sigma_2) \label{vs3}\ee
where
\be s_1(x)=\frac{f''(\xi_1(x))}{2} (x+1-\sigma_1).\ee
The first derivatives of \eqref{vs3} are
\begin{eqnarray} 
\nonumber v'(x)&=&s_1(x)+s_1'(x)(x-1-\sigma_2)\\
\nonumber v''(x)&=&2s_1'(x)+s_1''(x)(x-1-\sigma_2)\\
\nonumber \dots\\
\nonumber v^{(m)}(x)&=&ms_1^{(m-1)}(x)+s_1^{(m)}(x)(x-1-\sigma_2)
\end{eqnarray}
and evaluated in $x=1$:
\begin{eqnarray}
\nonumber v(1)&=&-s_1(1)\sigma_2\\
\nonumber v'(1)&=&s_1(1)-s_1'(1)\sigma_2\\
\dots\\
\nonumber v^{(m)}(1)&=&ms_1^{(m-1)}(1)-s_1^{(m)}(1)\sigma_2.
\end{eqnarray}
Similar results can be derived in the vicinity of the other endpoint.
We can then write \eqref{asyFilon} as
\begin{align*} 
I[v]\sim -\sum_{m=0}^\infty \frac{1}{(-{\rm i}\omega)^{m+1}}\Big[&e^{{\rm i}\omega}[ms_1^{(m-1)}(1)-s_1^{(m)}(1)\sigma_2]\\
&-e^{-{\rm i}\omega}[m{s}_{-1}^{(m-1)}(-1)-{s}_{-1}^{(m)}(-1)\sigma_1]\Big].
\end{align*}
If we collect the first terms in $s_1(1), s_1'(1), s_1''(1),...$, i.e.\ in $f^{''}(\xi_1(1))$, $f^{(iii)}(\xi_1(1))$, $f^{(iv)}(\xi_1(1)),\dots$, we have
\begin{align}
\nonumber I[v]\sim& s_1(1)e^{{\rm i}\omega}\left[\frac{\sigma_2}{{\rm i}\omega}-\frac{1}{\omega^2}\right]+s_1'(1)e^{{\rm i}\omega}\left[\frac{\sigma_2}{\omega^2}+\frac{2}{{\rm i}\omega^3}\right]+s_1''(1)e^{{\rm i}\omega}\left[
-\frac{\sigma_2}{{\rm i}\omega^3}+\frac{3}{\omega^4}
\right]+\dots\\\nonumber
&+{s}_{-1}(-1)e^{-{\rm i}\omega}\left[\frac{\sigma_1}{{\rm i}\omega}-\frac{1}{\omega^2}\right]+{s}_{-1}'(-1)e^{-{\rm i}\omega}\left[\frac{\sigma_1}{\omega^2}+\frac{2}{{\rm i}\omega^3}\right]\\&+{s}_{-1}''(-1)e^{-{\rm i}\omega}\left[
-\frac{\sigma_1}{{\rm i}\omega^3}+\frac{3}{\omega^4}
\right]+\dots
\label{sba}
\end{align}  
When we assume antisymmetric quadrature nodes ${c}_1=-{c}_2$, i.e.\ $\sigma_1=-\sigma_2$, the asymptotic order cannot be better than two: the coefficients of $s_1(1)$ and ${s}_{-1}(-1)$ in \eqref{sba} cannot both vanish in this case.
If we drop, however, the restriction ${c}_1=-{c}_2$, we have an extra free parameter and we might be able to obtain better asymptotic order for carefully chosen ${c}_1=-1+\sigma_1(\omega)$ and ${c}_2=1+\sigma_2(\omega)$ nodes. The equations
\be
\frac{\sigma_1}{{\rm i}\omega}-\frac{1}{\omega^2}=0,\quad \frac{\sigma_2}{{\rm i}\omega}-\frac{1}{\omega^2}=0
\ee
are satisfied for $\sigma_1=\sigma_2=i/\omega$. From \eqref{sba} it is clear that the choice ${c}_1=-1+i/\omega$ and ${c}_2=1+i/\omega$ gives us an $O(\omega^{-3})$ scheme for $\omega\gg 1$ and that asymptotic order 3 is the best asymptotic behaviour one may reach using a Filon-type method with two nodes. 
The quadrature nodes of this complex Filon-type method can also obtained using a steepest descent analysis, see \cite{Deano2009,Huybrechs2006,Huybrechs2010}.
 The quadrature rule resulting is 
\begin{equation}\label{complex}
Q_2^{F\text{-}C}=\frac{{\rm i}\left[e^{-{\rm i}\omega}f(-1+{\rm i}/\omega)-e^{{\rm i}\omega}f\left(1+{\rm i}/\omega\right)\right]}{\omega}.
\end{equation}
Results are shown in Figure \ref{fig:filon2}. Clearly the $Q_2^{F\text{-}C}$ method  shows the best behaviour for large $\omega$ values.  When using only real interpolation nodes, one can only construct methods with asymptotic order two. In this case, one can use interpolation points such as \eqref{adfilon} to construct an adaptive Filon method $Q_2^{F\text{-}A}$ as an alternative to the interpolation points of the EF approach.

\begin{figure}
	\centering
			\includegraphics[width=0.7\textwidth]{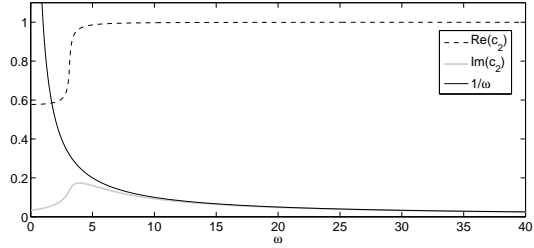}
			\caption{Real and imaginary part of $c_2(\omega)$, one of the quadrature nodes of the $Q_2^{F\text{-}CA}$ rule.}
	\label{fig:fca}
\end{figure}

The complex Filon method $Q_2^{F\text{-}C}$ should not be used for small $\omega$ values. As can be seen in Figure \ref{fig:filon2}, the $Q_2^{F\text{-}C}$ method does not reduce to the classical Gaussian scheme for $\omega\to 0$. One can, however, construct a Filon method  $Q_2^{F\text{-}CA}$ which combines the high asymptotic order of $Q_2^{F\text{-}C}$ with good results for smaller frequencies. Where the real parts of the interpolation nodes of the $Q_2^{F\text{-}C}$ method equal the endpoints over the whole $\omega$ range, we let now vary the real part of the interpolation nodes $c_1$ and $c_2$ from the Legendre nodes ${\bar c}_1$ and ${\bar c}_2$ to the endpoints $-1$ and $1$. For the imaginary part we take a curve which starts off close to zero for small $\omega$ but soon goes to the $i/\omega$ curve for larger $\omega$. We used the following $S$-functions (see Figure \ref{fig:fca}) to produce the $Q_2^{F\text{-}CA}$ results shown in Figure \ref{fig:filon2}:
\begin{equation}\label{ReF}
\nonumber {\rm Re}(c_2(\omega))=-{\rm Re}(c_1(\omega))=({{\bar c}_2-1)S(\omega;\pi;2)}+1,
\end{equation}
and
\begin{equation}\label{ImF}
\nonumber {\rm Im}(c_1(\omega))={\rm Im}(c_2(\omega))=\frac{1-S(\omega;\pi;1)}{\omega},
\end{equation}
with $S(\omega;r;n)$ as in Eq.\ \eqref{adfilon}.

\subsection{Methods with $\nu=3$ quadrature nodes}\label{sectionnu3}
\begin{figure}
	\centering
				\includegraphics[width=0.8\textwidth]{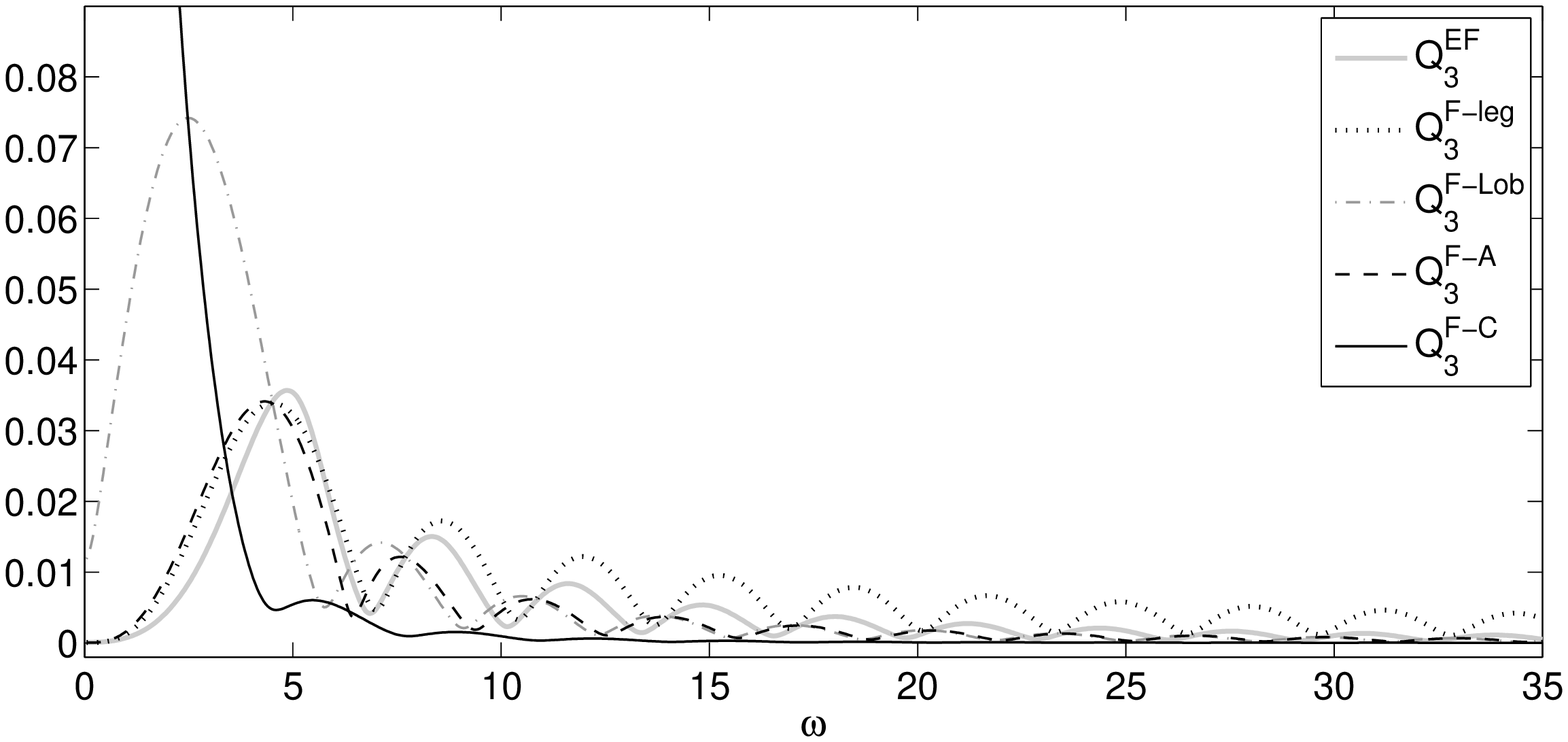}
					\includegraphics[width=0.8\textwidth]{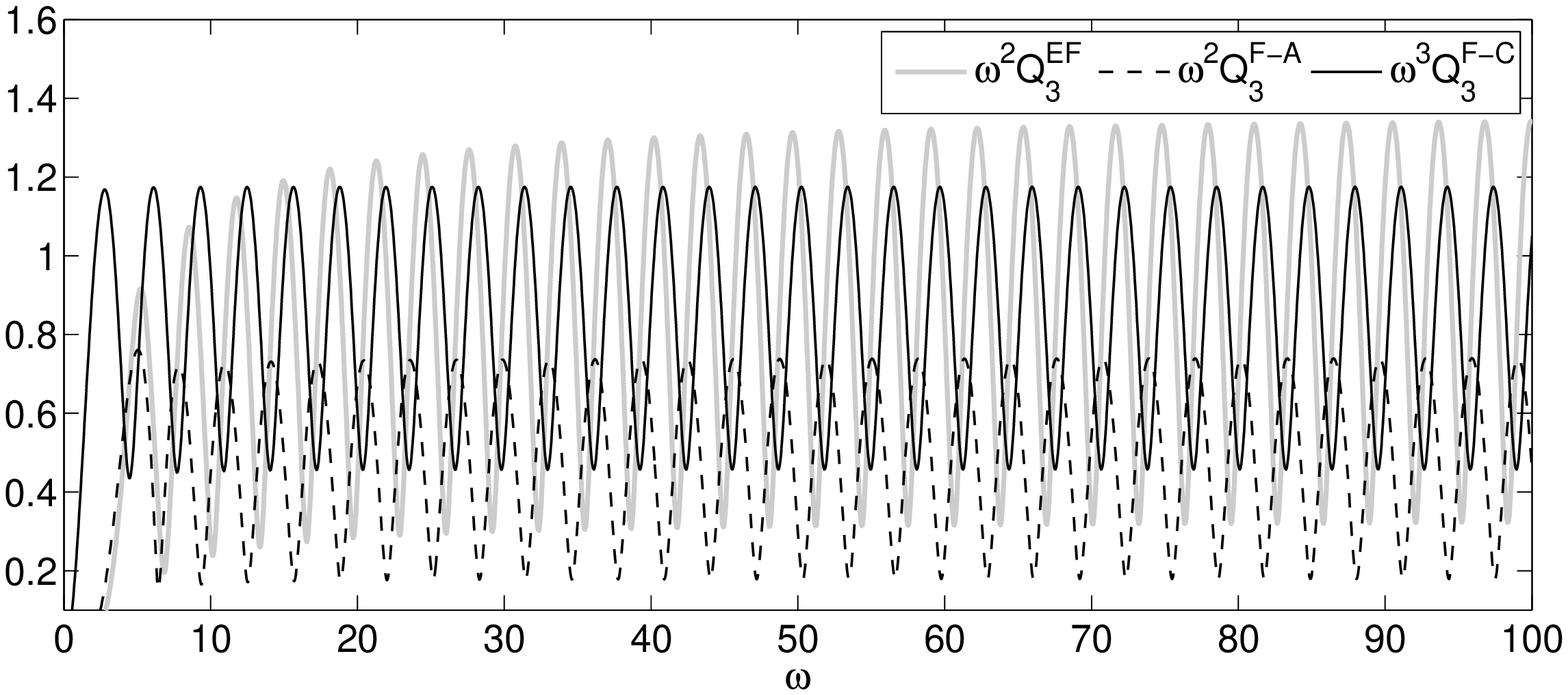}
	\caption{The error in some $\nu=3$ Filon-type schemes for $\int_{-1}^1 e^x e^{{\rm i}\omega x}dx$ and different values of $\omega$. The top graph shows the absolute error and the bottom figure contains the normalised errors.}\label{fig:fig3}
\end{figure}

The EF method has ${c}_1=-{c}_3$, ${c}_2$=0 and $w_1=w_3$. An iteration procedure is used to solve the system
\[\int_{-1}^1 x ^P e^{{\rm i}\omega x}dx-w_1 {c}_1^P e^{{\rm i}\omega {c}_1}-w_1 (-{c}_1)^P e^{-{\rm i}\omega {c}_1}-w_2 0^P=0,\quad P=0,1,2\]
The problems related to the ill-conditioning, mentioned before for the $\nu=2$ EF rule, continue to exist, making it difficult to compute accurate nodes and weights for very large $\omega$. The EF nodes and weights are shown in Figure \ref{fig:figef2b}. The EF scheme can be interpreted as a Filon method with the EF interpolation nodes $c_1$, $c_2$ and $c_3$. This Filon/EF scheme $Q_3^{EF}$ is compared with a Filon method with (fixed) Lobatto or Legendre nodes in Figure \ref{fig:fig3}. 
As for the $\nu=2$ Filon method, it can be shown that no asymptotic order higher than two can be reached for a Filon-type method with $c_1=-c_3$, $c_2=0$. 
One can again avoid the (expensive and ill-conditioned) EF-iteration procedure by using alternative interpolation nodes to construct an adaptive Filon rule $Q_3^{F\text{-}A}$. We use the formula \eqref{adfilon} but now with the Legendre node ${\bar c_3}=\sqrt{3/5}$ (see Figure \ref{fig:figc1}). 

As for the $\nu=2$ Filon rule, the asymptotic order can be increased by allowing complex interpolation nodes. We use the same interpolation nodes as for the $\nu=2$ case around the endpoints, i.e.\ $c_1(\omega)=-1+{\rm i}/\omega,  c_3(\omega)=1+{\rm i}/\omega$. The value of $c_2$ has no influence on the asymptotic order of the method and we chose to take the middle of the interval: $c_2=0$. Figure \ref{fig:fig3} shows numerical results for both the adaptive Filon rule $Q_3^{F\text{-}A}$ and the rule $Q_3^{F\text{-}C}$ with complex nodes. Going from two to three interpolation nodes, does not increase the asymptotic order (the asymptotic order is determined by the number of interpolation nodes in the vicinity of the endpoints) but does make the magnitude of the error smaller by two decimal units. Another natural choice for the complex Filon rule may be to take the interpolation node $c_2={\rm i}/\omega$. This choice leads to results very similar to the ones for the $Q_3^{F\text{-}C}$ rule shown in Figure \ref{fig:fig3}: for larger frequencies $\omega$ the error behaviours are identical, only for small $\omega < 2\pi$ --where other rules should be prefered-- there is some difference.

\subsection{Methods with $\nu=4$ quadrature nodes}\label{nu4}
Figure \ref{fig:figef2b} shows the nodes ${c}_1=-{c_4},{c}_2=-{c_3}$ and weights $w_1=w_4, w_2= w_3$ of the EF version of the four-node rule. 
\begin{figure}
	\centering
		\includegraphics[width=0.7\textwidth]{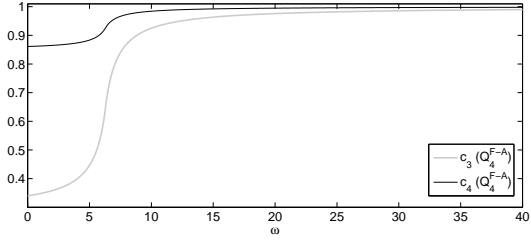}
	\caption{Quadrature nodes $c_3(\omega)$ and $c_4(\omega)$ of the adaptive Filon method $Q_4^{F\text{-}A}$.}
	\label{fig:nodes4}
\end{figure}
As for the $n=2,3$ schemes, an alternative to this EF rule exists in the form of an adaptive Filon method $Q_4^{F\text{-}A}$. The interpolation nodes of this adaptive Filon method are shown in Figure \ref{fig:nodes4}. Again the $S$-shaped function \eqref{adfilon} is used. This method has asymptotic order three. Asymptotic order five can be reached when using the complex interpolation nodes
\begin{eqnarray}
&c_1(\omega)=-1+(2-\sqrt{2}){\rm i}/\omega, &\;\;\;\;c_2(\omega)=-1+(2+\sqrt{2}){\rm i}/\omega,\label{complexnodes0}\\
&c_3(\omega)=1+(2-\sqrt{2}){\rm i}/\omega, &\;\;\;\;c_4(\omega)=1+(2+\sqrt{2}){\rm i}/\omega,
\label{complexnodes}
\end{eqnarray}
to construct the $Q_4^{F\text{-}C}$ method. Higher asymptotic orders are not possible for a method with four interpolation nodes. This can again be shown using the asymptotic expression for the error of the Filon method \eqref{asyFilon}. We have now two interpolation nodes in the vicinity of $-1$, which we denote by $-1+\sigma_1(\omega)$ and $-1+\sigma_2(\omega)$ and two interpolation nodes in the vicinity of $1$ denoted by $1+\sigma_3(\omega)$ and $1+\sigma_4(\omega)$. The interpolation error in the vicinity of $1$ is then
\[v(x)=s_1(x)(x-1-\sigma_3)(x-1-\sigma_4)\]
with
\[s_1(x)=\frac{f^{(iv)}(\xi_1(x))}{4!}(x+1-\sigma_1)(x+1-\sigma_2).\]
The first derivatives of $v(x)$ are
\begin{eqnarray}
\nonumber v'(x)&=&s_1(x)S(x)+s_1'(x)P(x)\\
\nonumber v''(x)&=&2s_1(x)+2s_1'(x)S(x)+s_1''(x)P(x)\\
\nonumber v^{(iii)}(x)&=&6s_1'(x)+3s_1''(x)S(x)+s_1^{(iii)}(x)P(x)\\
\nonumber v^{(iv)}(x)&=&12s_1''(x)+4s_1^{(iii)}(x)S(x)+s_1^{(iv)}(x)P(x)\\
&\dots
\end{eqnarray}
with $S(x)=(x-1-\sigma_3)+(x-1-\sigma_4)$ and $P(x)=(x-1-\sigma_3)(x-1-\sigma_4)$.
This gives us
\begin{eqnarray}
\nonumber v(1)&=&s_1(1)\sigma_3\sigma_4\\
\nonumber v'(1)&=&-s_1(1)[\sigma_3+\sigma_4]+s_1'(1)\sigma_3\sigma_4\\
\nonumber v''(1)&=&2s_1(1)-2s_1'(1)[\sigma_3+\sigma_4]+s_1''(1)\sigma_3\sigma_4\\
\nonumber v^{(iii)}(1)&=&6s_1'(1)-3s_1''(1)[\sigma_3+\sigma_4]+s_1^{(iii)}(1)\sigma_3\sigma_4\\
\nonumber v^{(iv)}(1)&=&12s_1''(1)-4s_1^{(iii)}(1)[\sigma_3+\sigma_4]+s_1^{(iv)}(1)\sigma_3\sigma_4\\
&\dots
\end{eqnarray}
The asymptotic expression \eqref{asyFilon} for the Filon error then reduces to the following form:
\begin{align}\label{asymptexp4}
I[v]\sim& \nonumber s_1(1)e^{{\rm i}\omega}\left[-\frac{\sigma_3\sigma_4}{{\rm i}\omega}+\frac{\sigma_3+\sigma_4}{\omega^2}+\frac{2}{{\rm i}\omega^3}\right]-s_1'(1)e^{{\rm i}\omega}\left[\frac{\sigma_3\sigma_4}{\omega^2}+\frac{2(\sigma_3+\sigma_4)}{{\rm i}\omega^3}-\frac{6}{\omega^4}\right]\\\nonumber
&+s_1''(1)e^{{\rm i}\omega}\left[\frac{\sigma_3\sigma_4}{{\rm i}\omega^3}-\frac{3(\sigma_3+\sigma_4)}{\omega^4}-\frac{12}{{\rm i}\omega^5}\right]+\dots\\
\nonumber&+{s}_{-1}(-1)e^{-{\rm i}\omega}\left[-\frac{\sigma_1\sigma_2}{{\rm i}\omega}+\frac{\sigma_1+\sigma_2}{\omega^2}+\frac{2}{{\rm i}\omega^3}\right]\\\nonumber&-{ s}_ {-1}'(-1)e^{-{\rm i}\omega}\left[\frac{\sigma_1\sigma_2}{\omega^2}+\frac{2(\sigma_1+\sigma_2)}{{\rm i}\omega^3}-\frac{6}{\omega^4}\right]\\
&+{s}_{-1}''(-1)e^{-{\rm i}\omega}\left[\frac{\sigma_1\sigma_2}{{\rm i}\omega^3}-\frac{3(\sigma_1+\sigma_2)}{\omega^4}-\frac{12}{{\rm i}\omega^5}\right]+\dots
\end{align}
To have an $O(\omega^{-5})$ scheme, the following equations need to be satisfied (and similarly for $\sigma_3$ and $\sigma_4$):
\begin{eqnarray}
-\frac{\sigma_1\sigma_2}{{\rm i}\omega}+\frac{\sigma_1+\sigma_2}{\omega^2}+\frac{2}{{\rm i}\omega^3}=0 \label{s31}\\
\frac{\sigma_1\sigma_2}{\omega^2}+\frac{2(\sigma_1+\sigma_2)}{{\rm i}\omega^3}-\frac{6}{\omega^4}=0 \label{s32}
\end{eqnarray}
This gives us $\sigma_1(=\sigma_3)=(2-\sqrt{2}){\rm i}/\omega$ and $\sigma_2(=\sigma_4)=(2+\sqrt{2}){\rm i}/\omega$ and an $O(\omega^{-5})$ asymptotic error. This corresponds again to the nodes given by the numerical steepest descent approach from \cite{Huybrechs2006}. Indeed, combining \eqref{s31} and \eqref{s32} gives us  $-\omega^2 \sigma_2+4{\rm i}\omega \sigma_2+2= L_2(\omega \sigma_2/{\rm i})=0$ where $L_2$ is the Laguerre polynomial of degree two.

\begin{figure}
	\centering
				\includegraphics[width=0.8\textwidth]{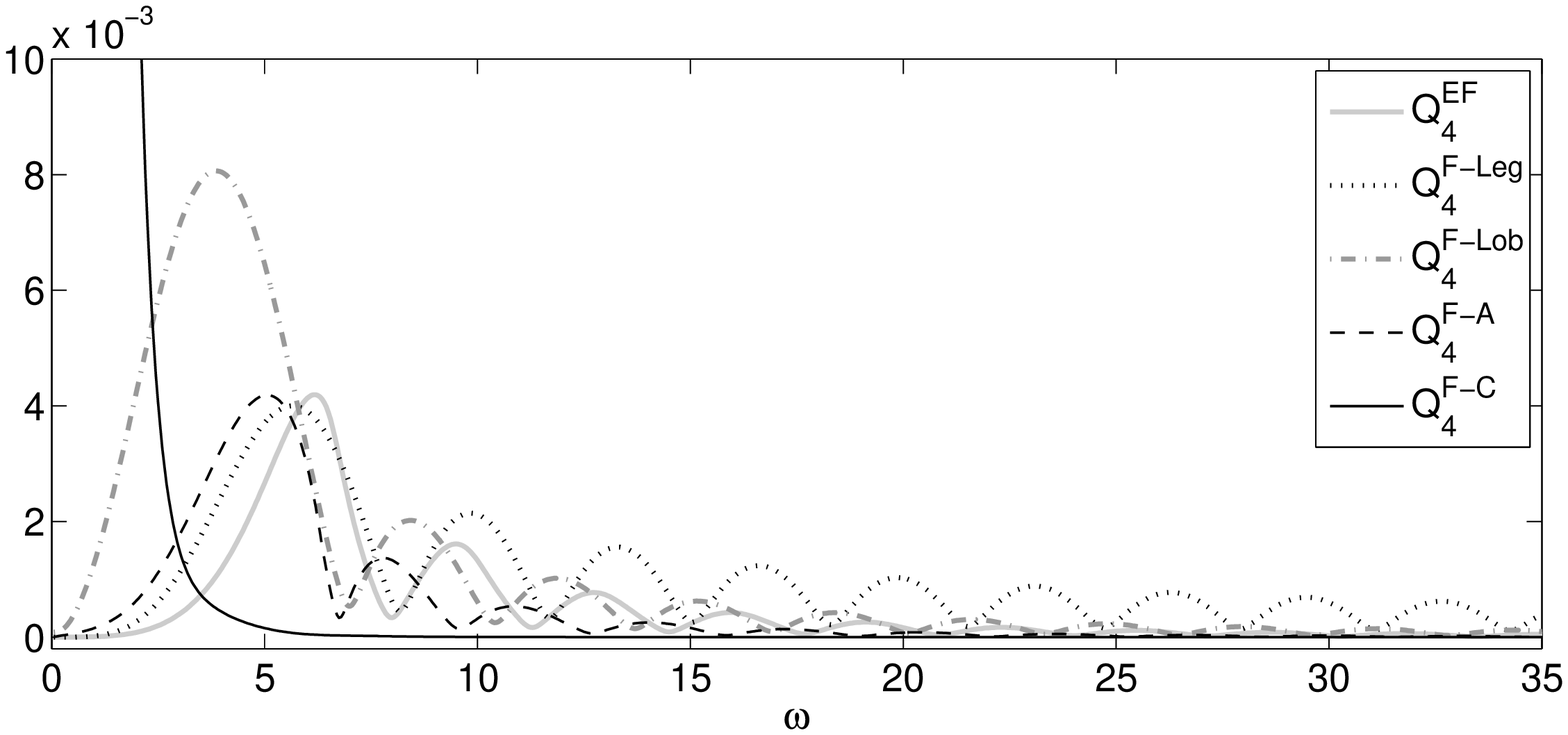}
								\includegraphics[width=0.8\textwidth]{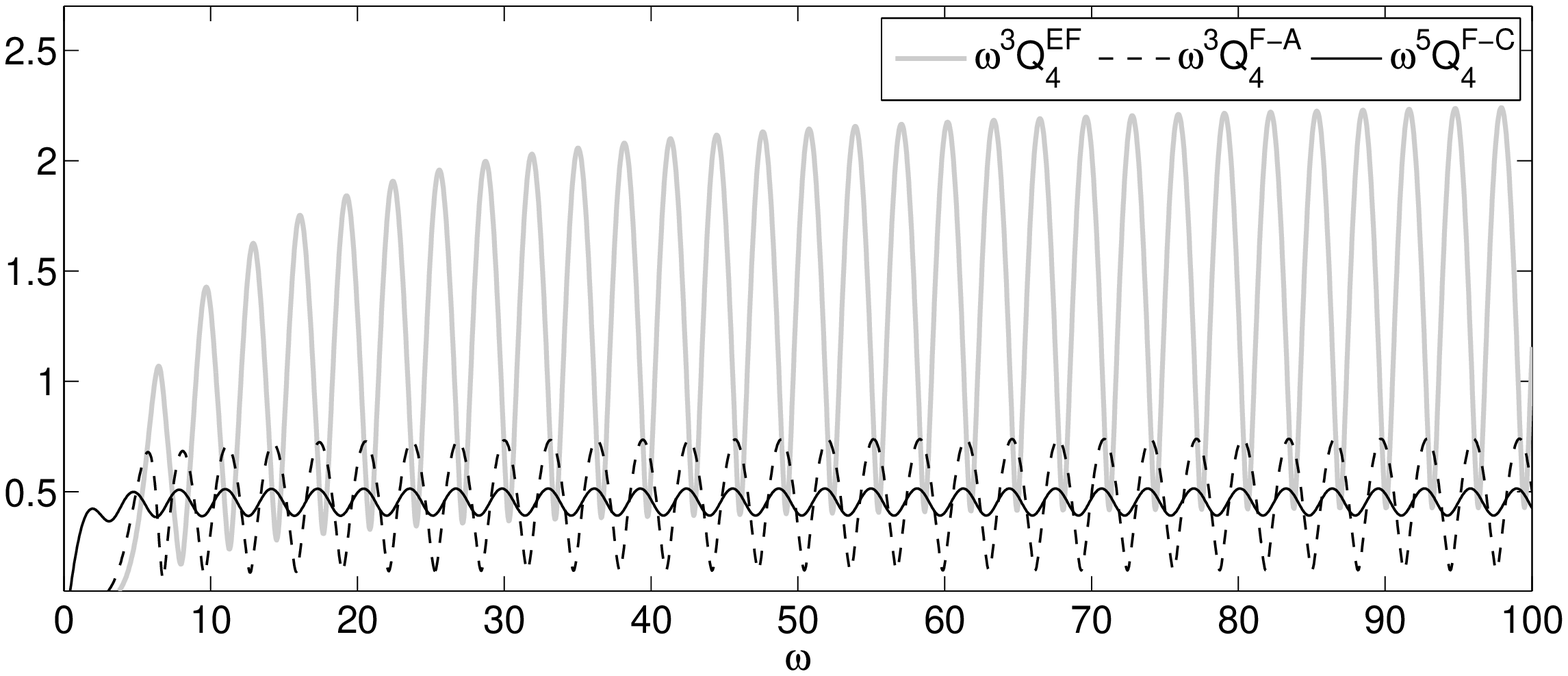}
	\caption{The error in some $\nu=4$ Filon-type schemes for $\int_{-1}^1 e^x e^{{\rm i}\omega x}dx$ and different values of $\omega$. The top graph shows the absolute error and the bottom figure contains the normalised errors.}\label{fig:figfilon4}
\end{figure}

Figure \ref{fig:figfilon4} shows the results of the different four node schemes applied on our test problem.

\subsection{Methods with more quadrature nodes}
\begin{figure}
	\centering
		\includegraphics[width=0.8\textwidth]{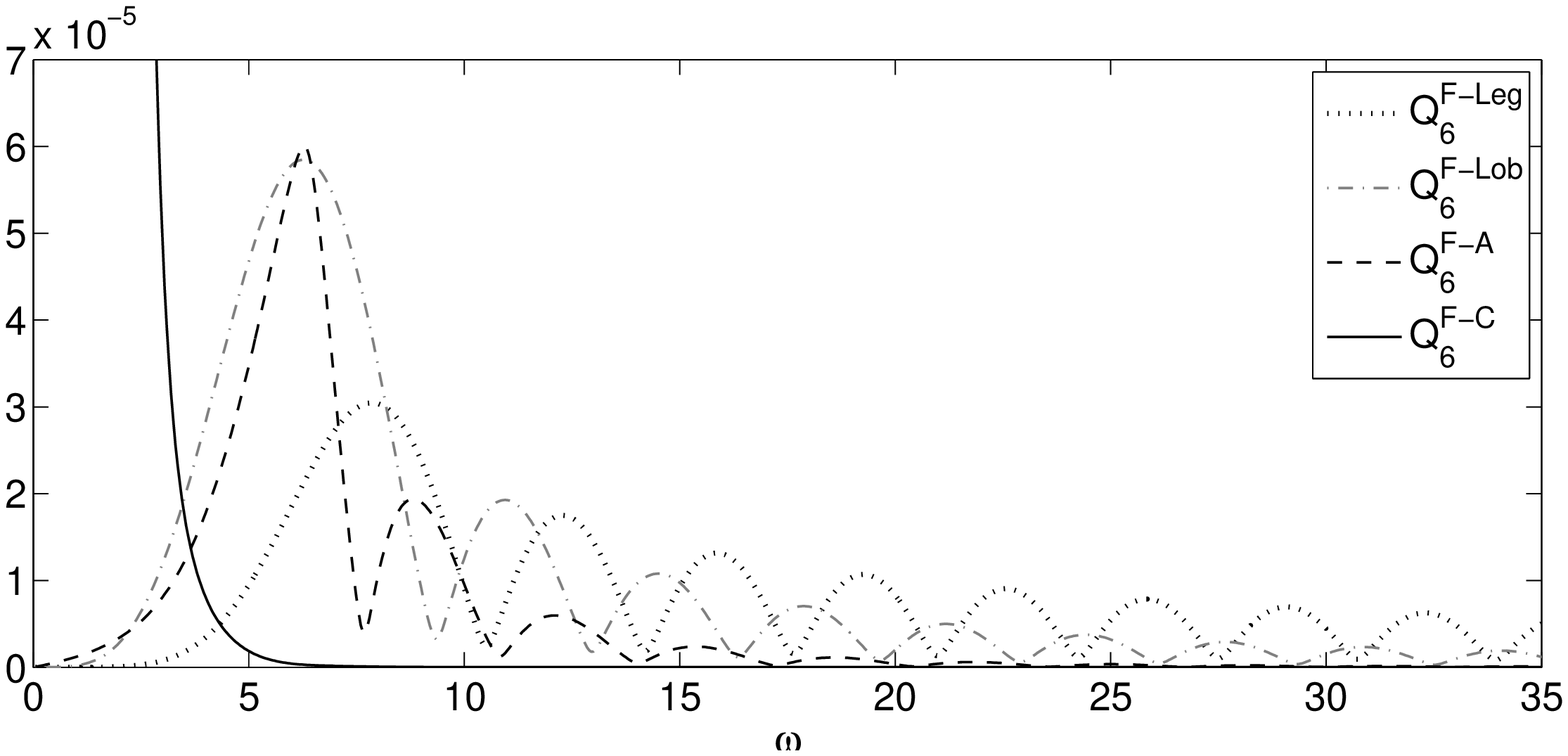}
				\includegraphics[width=0.8\textwidth]{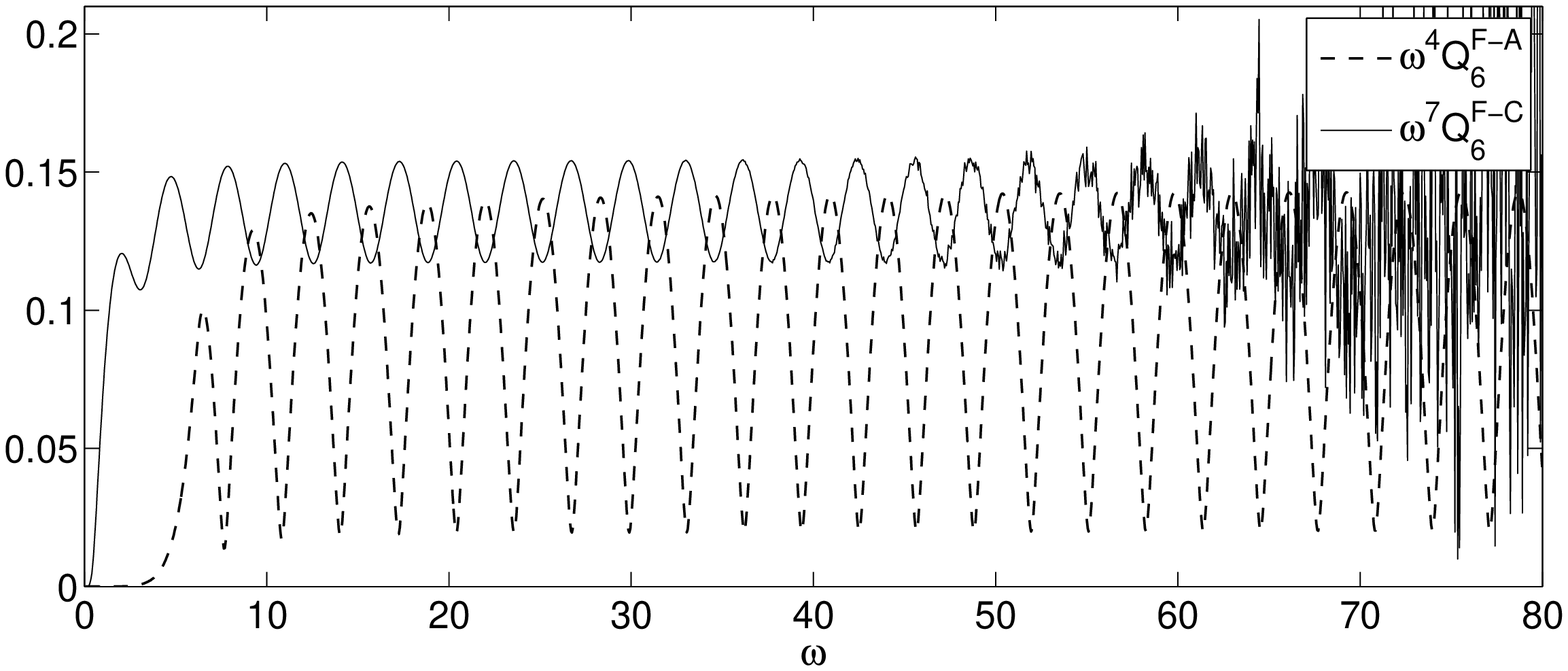}
	\caption{The error in some $\nu=6$ Filon-type schemes for $\int_{-1}^1 e^x e^{{\rm i}\omega x}dx$ and different values of $\omega$.
	}\label{fig:figfilon6}
\end{figure}

The asymptotic order can be increased further by introducing more quadrature points. Figure \ref{fig:figfilon6} shows the results for some methods with $\nu=6$ quadrature/interpolation nodes. 
The Filon-type schemes, as used in figure \ref{fig:figfilon6}, are constructed by replacing $f$ in $I[f]=\int_{-1}^1 f(x)e^{i\omega x}dx$ by a Lagrange interpolating polynomial $p=\sum_{j=0}^{\nu-1} a_j x^j$ so that 
\be
I[f]\approx I[p]=\sum_{j=0}^{\nu-1} a_j I[x^j]
\label{ac}\ee where the terms $I[x^j]$ are explicitly given by \eqref{Gamma}. The coefficients $a_j$ of the polynomial $p$ depend on the position of the quadrature nodes. The coefficients $a_0$ and $a_{\nu-1}$, for instance, are given by
\be
a_0=\sum_{j=1}^{\nu}\frac{f(c_j)\prod_{k=1,k\neq j}^{\nu}c_k}{\prod_{k=1,k\neq j}^{\nu}(c_j-c_k)}\label{a0}
\ee
and 
\be
a_{\nu-1}=\sum_{j=1}^{\nu}\frac{f(c_j)}{\prod_{k=1,k\neq j}^{\nu}(c_j-c_k)}.\label{av}
\ee
Unfortunaly however, more and more numerical difficulties due to cancellation effects appear in this procedure when $\nu$ is increased, especially for large values of $\omega$ where the nodes are closely clustered near the endpoints, in which case the coefficients $a_j$ in \eqref{ac}, e.g.\ \eqref{a0} and \eqref{av}, have large values. This means that in a practical implementation there is a limit on the highest asymptotic order which is feasible. In the next section, it will be shown, however, how the accuracy of our methods can be further increased by adding extra interior nodes.

\section{An automatic quadrature scheme}
When applying a quadrature rule into an automatic software package in which an integral needs to be approximated within a prescribed tolerance, one traditionally uses a process in which the integral is approximated using static quadrature rules on adaptively refined subintervals of the integration domain. This adaptive quadrature technique is very effective for non-oscillatory integrands (small values of $\omega$).
\begin{figure}
	\centering
	\includegraphics[width=0.8\textwidth]{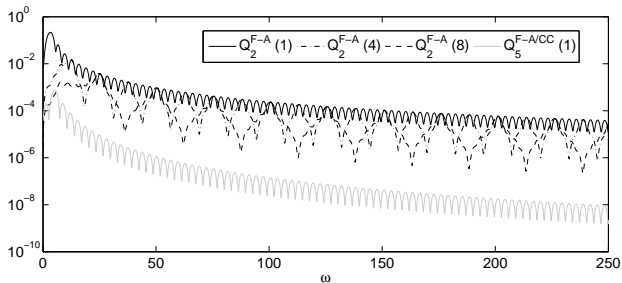}
	\caption{The error in some Filon-type schemes applied on $\int_{-1}^1 e^xe^{i\omega x}dx$. The notation $Q_2^{F\text{-}A}(N)$ is used to denote the $Q_2^{F\text{-}A}$ method applied on $N$ equidistant subintervals of the interval $[-1,1]$. The $Q_5^{F\text{-}A/CC}(1)$ method interpolates through 5 nodes: the 2 nodes from $Q_2^{F\text{-}A}(1)$ and the `Chebyshev' nodes $-1$, $0$ and $1$. }\label{fig:nietSplitsen}
\end{figure}
In the oscillatory regime, however, this classical way of thinking does not lead to good results at all.
This is illustrated in Figure \ref{fig:nietSplitsen}. The problem $\int_{-1}^1 e^xe^{i\omega x}dx$ is solved by the $Q_2^{F\text{-}A}$ method which is applied ($i$) once on the whole interval (2 function evaluations), ($ii$) on each of 4 subintervals of equal length (8 function evaluations), ($iii$) on each of 8 subintervals of equal length (16 function evaluations). The results are compared with the $Q_5^{F\text{-}A/CC}$ method which is a Filon method interpolating $f$ through the 2 nodes of the  $Q_2^{F\text{-}A}$ rule and the 3 extra nodes $-1$, $0$, and $1$. The $Q_5^{F\text{-}A/CC}$ scheme is applied once on the whole integration interval and needs (only) 5 evaluations of the function $f$. The three different applications of the $Q_2^{F\text{-}A}$ method display a $O(\omega^{-2})$ behaviour, while $Q_5^{F\text{-}A/CC}$ has a $O(\omega^{-3})$ asymptotic error. But a more important observation is that subdividing the integration interval and applying the $Q_2^{F\text{-}A}$ quadrature rule on each subinterval does not lead to a decrease in the maximal value of the error envelope. Intuitively this is not a surprise considering the asymptotic expansion \eqref{asy}. The value of an oscillating integral is determined mostly by the regions where the oscillations do not cancel out \cite{Cools}, i.e.\ the most important contribution to the integral comes from the endpoints. Dividing the integration interval into subintervals introduces extra (interior) endpoints, which means that one has now to include also sufficient information around these interior endpoints to obtain a good accuracy in the approximation of each subintegral. 
For oscillatory integrals, clearly, better results are obtained by introducing extra nodes into the quadrature rule rather than subdividing the integration interval. One can for instance combine our $\omega$ dependent nodes with the $n$ Chebyshev nodes (in $[-1,1]$)
\[x_j=\cos\frac{j\pi}{n-1},\;\;0\leq j\leq n-1.\]
Including these Chebyshev nodes, one also controls the interpolation error and ensures the convergence of the scheme (see \cite{Evans1999,Melenk}). The $\omega$ dependent nodes improve the asymptotic behaviour and the combination with the Chebyshev nodes has some extra advantages (see \cite{Dominguez}): ($i$) it allows to express the polynomial approximation of $f$ as a Chebyshev expansion avoiding the cancellation issues which appear in the Lagrange polynomial interpolation for larger numbers of nodes ($ii$) another particular attraction of the Chebyshev points is the fact that they are nested and we have thus a convenient reuse of nodes when the number of nodes is doubled ($iii$) the quadrature weights corresponding to $n$ Chebyshev points can be evaluated in $O(n\log n)$ time by FFT algorithms \cite{Trefethen2008}.
A so-called Filon-Clenshaw-Curtis rule, based on replacing $f$ by a polynomial interpolant through the Chebyshev points, was introduced in \cite{Dominguez}.  The technique presented in this section will in fact improve the asymptotic behaviour of these Filon-Clenshaw-Curtis rules, which means that a smaller number of iterations (and consequently a smaller number of function evaluations) will be needed to reach a prescribed tolerance. Similar ideas were also discussed in \cite{Huybrechs2010}. In \cite{Huybrechs2010}, some so-called superinterpolation nodes ensuring optimal asymptotic behaviour are combined with the nodes used by Fej\'er in his first quadrature rule to form a complex Filon method. The nodes are, however, not nested and thus less suited for an automatic quadrature scheme. We adapt here the Algorithm described in \cite{Huybrechs2010} to the inclusion of the Chebyshev points, giving us Algorithm \ref{Algo1} which forms the Filon-type method interpolating $f$ at the $n$ Chebyshev nodes, i.e. the $n$ extreme points of the Chebyshev polynomial 
$T_{n-1}(x)$ in $[-1,1]$, including the endpoints, and $\nu$ $\omega$-dependent nodes $c_1,\dots,c_\nu$. 

\noindent\begin{myAlgo}\label{Algo1}{ {\rm Adaptive Filon-type method}}\\
Input: f, integer $n$ and $\nu$ additional points ${\bf c}=(c_1,\dots,c_\nu)^T$.\vspace*{1mm}
\begin{enumerate}
\item Interpolate $f$ by a polynomial $p_1={\bf T_{n-1}}{\bf p_1}$ through $n$ Chebyshev points using the discrete cosine transform (see \cite{Dominguez}), which, by the FFT, requires $O(n\log n)$ operations.
${\bf T_{n-1}}$ is defined as $(T_0,T_1,\dots,T_{n-1})$ and ${\bf p_1}=(p_{1,1},\dots,p_{1,n})^T$.
\vspace*{3mm}
\item Interpolate $\displaystyle g=\frac{f-p_1}{(x^2-1)T'_{n-1}}=\frac{2(f-p_1)}{(n-1)(T_n-T_{n-2})}$ through the points ${\bf c}$ by a polynomial $p_2={\bf T_{\nu-1}}{\bf p_2}$.
\vspace*{3mm}
\item Form $p_3={(x^2-1)T'_{n-1}}p_2=(n-1)(T_n-T_{n-2})p_2/2={\bf T_{n+\nu-1}}{\bf p_3}$ using the identity
\[2T_rT_s=T_{r+s}+T_{|r-s|}.\]
Note that $p_3$ vanishes at the $n$ Chebyshev points and equals $f-p_1$ for the nodes ${\bf c}$.
\vspace*{2mm}
\item The polynomial $p$ which interpolates through $n$ Chebyshev points and the nodes ${\bf c}$ is then
\[p=p_1+p_3={\bf T_{n+\nu-1}}[{\bf p_1}+{\bf p_3}]\]
where ${\bf p_1}$ is padded at the end with zeros. The adaptive Filon method is then defined as
\[I[p]=\sum_{k=0}^{n+\nu-1}{\bf\tau}_{k}[p_{1,k+1}+p_{3,k+1}]\]
where
\[\tau_{k}:\int_{-1}^1T_k(x)e^{i\omega x}dx.\]

\end{enumerate}
\vspace*{3mm}
\end{myAlgo}
A robust implementation of the method requires an efficient way of evaluating the `moments' $\tau_{k}$. A recursive algorithm for computing $\tau_{k}$ may be based on the recurrence relations for Chebyshev polynomials. However, this algorithm is only stable when $k\leq\omega$ and is inappropriate for evaluating $\tau_{k}$ when $k>\omega$. This problem has recently been adressed and resolved in \cite{Dominguez}. The values $\tau_k$ with $k>\omega$ are computed by adding a second phase to the algorithm in which a tridiagonal system of equations is solved. The right-hand side of this tridiagonal system is determined by an asymptotic argument and the resulting `two-phase' algorithm is accurate and stable for all $k$ and $\omega$. 

We propose two values for $\nu$: $\nu=2$ and $\nu=4$. Larger values for $\nu$ can be chosen but may be feasible only with high-precision calculations for very large $\omega$ values since one may have  numerical stability issues: significant cancellation appears when the nodes are very closely clustered near the endpoints. A good choice for the $\nu=2$ nodes $c_1$ and $c_2$ are the $S$-shaped nodes from eq.\ \eqref{adfilon} to form a method with asymptotic order 3 when combined with the Chebyshev nodes. As seen earlier, the asymptotic order is determined by the number of interpolation nodes in the vicinity of the endpoints. Here we have two interpolation nodes near each endpoint, i.e.\ the endpoint itself which is one of the Chebyshev nodes, and the $Q_2^{F-A}$ node \eqref{adfilon}. Consequently, the asymptotic order can be derived from eq.\ \eqref{asymptexp4}, where we take $\sigma_1=\sigma_4=0$ since the endpoints are now included among the interpolation points. From \eqref{asymptexp4} we also learn that asymptotic order 4 can be reached when allowing complex nodes which behave as $c_1=-1-2/(i\omega)$, $c_2=1-2/(i\omega)$ in the asymptotic regime.
In the case of $\nu=2$ frequency dependent nodes, the polynomials $p_2$ and $p_3$ appearing in Algorithm \ref{Algo1} have the following form:
\[p_2=p_{2,1}T_0+p_{2,2}T_1,\;{\rm with}\; p_{2,1}=\frac{c_1g(c_2)-c_2g(c_1)}{c_1-c_2},\; p_{2,2}=\frac{g(c_1)-g(c_2)}{c_1-c-2}\]
and $p_3={\bf T_{n+1}}{\bf p_3}$ with ${\bf p_3}=(p_{3,1},\dots,p_{3,n+2})^T$ and
\begin{eqnarray*}
   p_{3,n+2}=-p_{3,n-2}=\frac{(n-1)}{4}p_{2,2},\\
    p_{3,n+1}=-p_{3,n-1}=\frac{(n-1)}{2}p_{2,1},\\
    p_{3,n}=0,\;p_{3,j}=0\, (j=1,...,n-3).\end{eqnarray*}
For $\nu=4$ nodes, one can use the $S$-shaped $\omega$ dependent nodes from the $Q_4^{F-A}$ method, leading to asymptotic order 4. In this case, the coefficients of the polynomial $p_2$ are given by
\begin{eqnarray*}
  p_{2,1}&=&p_{2,3}-\sum_{m=1}^{4}\frac{(\prod_{r=1,r\neq m}^4 c_r) f(c_m)}{\sigma_m}\\
  p_{2,2}&=&3p_{2,4}+\sum_{m=1}^{4}\frac{\sum_{r=1,r\neq m}^4 c_r \sum_{s>r,s\neq m}^4c_s f(c_m)}{\sigma_m}\\
  p_{2,3}&=&-\frac{1}{2}\sum_{m=1}^{4}\frac{\sum_{r=1,r\neq m}^4 c_r f(c_m)}{\sigma_m}\\
  p_{2,4}&=&\frac{1}{4}\sum_{m=1}^{4}\frac{f(c_m)}{\sigma_m}\\
\end{eqnarray*}
with $\sigma_m=\prod_{s=1,s\neq m}^{4}(c_m-c_s)$ and the polynomial $p_3$ can be obtained as follows
\[p_3=\frac{(n-1)}{4}\sum_{k=0}^{\nu-1}p_{2,k+1}(-T_{n-2+k}-T_{|k-n+2|}+T_{n+k}+T_{|k-n|}).\]

\begin{figure}
	\centering
		\includegraphics[width=0.8\textwidth]{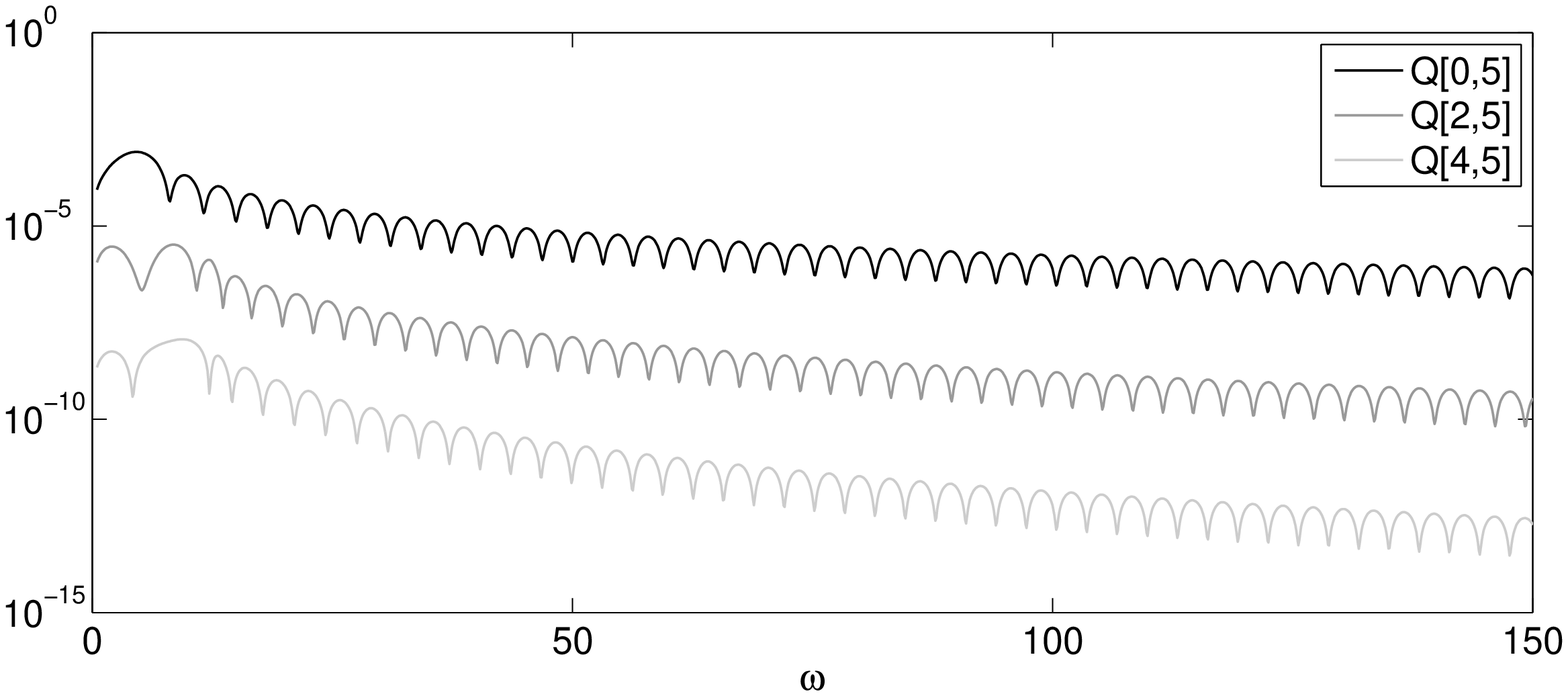}
		\includegraphics[width=0.8\textwidth]{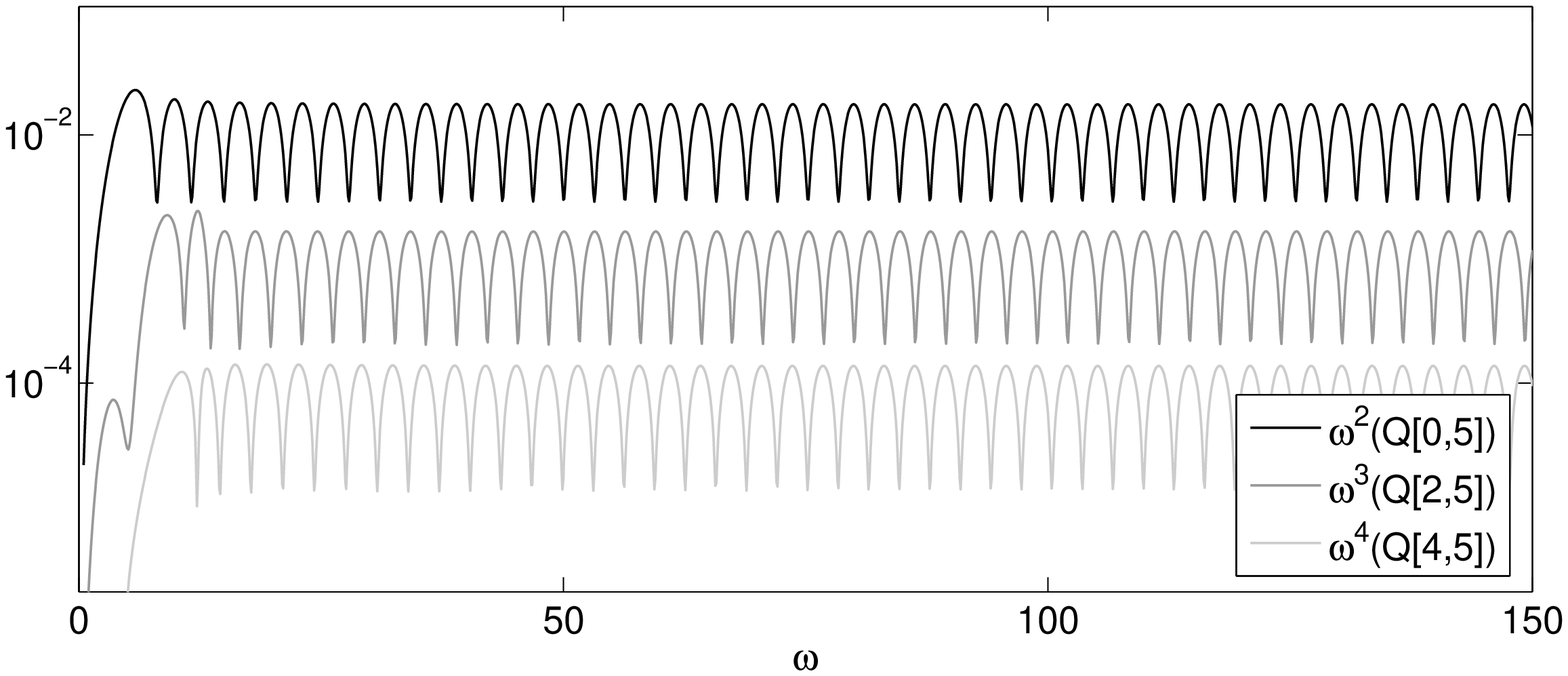}
	\caption{The error in some (Chebyshev) Filon-type schemes $Q[\nu,n]$ applied on $\int_{-1}^1 e^xe^{i\omega x}dx$ for different values of $\omega$. The top graph shows the absolute error and the bottom figure contains the normalised errors.}\label{fig:aut}
\end{figure}

Some numerical results are shown in Figure \ref{fig:aut}.
The notation $Q[0,5]$ is used to denote the Filon-Clenshaw-Curtis method \cite{Dominguez} with the 5 Chebyshev points $\cos(j\pi)/4$, $j=0,\dots,4$. The rules $Q[2,5]$ and $Q[4,5]$ implement Algorithm \ref{Algo1} and have $\nu=2$, resp. $\nu=4$ additional $\omega$ dependent nodes $c_j$, $j=1,\dots,\nu$. The $c_j$ nodes used in this figure, correspond to the $Q_2^{F-A}$ and the $Q_4^{F-A}$ nodes.

The following algorithm shows how the Filon method from Algorithm \ref{Algo1} can be applied in an adaptive procedure such that the estimated error in the eventually returned result is less than some tolerance specified by the user.

\begin{myAlgo}\label{Algo2}{\rm An automatic quadrature scheme}\\
Input: user tolerance $tol$, $f$, $\omega$.\\
Initial values: $n=3$ and ($\nu=2$ or $\nu=4$) nodes ${\bf c}$.
\begin{enumerate}
\item Form a first approximation $A_1$ using Algorithm 1 with the `asymptotic' nodes ${\bf c}$ and $n$ Chebyshev points.
\vspace*{1mm}
\item Form a second approximation $A_2$ using Algorithm 1 with the nodes ${\bf c}$ and $2n-1$ Chebyshev points. Only $n-1$ new evaluations of $f$ are needed.
\vspace*{1mm}
\item  If the difference between the two approximations $|A_1-A_2|$ is less than $tol$, accept the value returned by the first rule. Otherwise repeat steps 2-3 with $n:=2n-1$ and $A_1:=A_2$.
\end{enumerate}
\end{myAlgo}

In Table \ref{tab:tt} some results are listed for different values of $\omega$ for the problem with $f(x)=e^x$ on the integration interval $[-5,5]$. One can observe that the number of iterations, and consequently the number of function evaluations, decreases with $\omega$. Another observation is that the $\nu=2$ or $\nu=4$ rules clearly need a smaller number of function evaluations for larger $\omega$ than the $\nu=0$ method from \cite{Dominguez}: a consequence of the higher asymptotic order. The matlab implementation which was used to obtain these results is available at \verb| http://users.ugent.be/~vledoux/HOIsoftware|.

\begin{table}
	\caption{The true error, the number of iterations $n_{it}$ and the number of evaluations $n_{fev}$ of the function $f$ when applying Algorithm \ref{Algo2}, with $tol=10^{-9}$, on the problem $\int_{-5}^5 e^xe^{i\omega x}dx$. }
	\label{tab:tt}
		\begin{tabular}{cccc|ccc|ccc}
		\hline\noalign{\smallskip}
	&\multicolumn{3}{c}{$\nu=0$}&\multicolumn{3}{c}{$\nu=2$}&\multicolumn{3}{c}{$\nu=4$}\\
			$\omega$&$error$&$n_{it}$&$n_{fev}$&$error$&$n_{it}$&$n_{fev}$&$error$&$n_{it}$&$n_{fev}$\\
 \noalign{\smallskip}\hline\noalign{\smallskip}
 		10&1.5E-14&5&65&1.8E-10&4&35&3.8E-12&4&37\\
 		$100$&9.5E-11&4&33&7.8E-14&4&35&6.8E-11&3&21\\
 		$500$&3.2E-12&4&33&4.4E-10&3&19&7.5E-10&1&9\\
 		$1000$&1.2E-12&4&33&6.2E-11&3&19&4.3E-11&1&9\\
 		$5000$&3.7E-14&4&33&4.6E-10&1&7&7.5E-14&1&9\\
\noalign{\smallskip}\hline
		\end{tabular}
\end{table}

\section{Other remarks}
In literature \cite{Iserles2004b,Iserles2005}, Filon-type methods have been presented for the more general problem
\be\int_{-1}^1f(x)e^{{\rm i}\omega g(x)}dx\label{gen}\ee
An important distinction has then to be made between oscillators $g$ such that $g'\neq 0$ in $[-1,1]$ and oscillators with so-called stationary points $\xi\in[-1,1]$ where the $g'$ vanishes.
Once $g'$ vanishes at one or more points in the integration interval, optimal asymptotic behaviour can only be reached by taking into account not only information about $f$ in the endpoints but also in the stationary points (see \cite{Iserles2004b,Huybrechs2009}). A generalization of the methods from the previous sections consists, consequently, in the addition of some extra $\omega$ dependent nodes which tend to the stationary point as $\omega\to\infty$. Also these nodes can follow an $S$-shaped curve as they evolve from small to large $\omega$. The number of nodes to be added depends on the particular asymptotic order required and on the degree of the stationary point \cite{Iserles2004b}. 

The construction of a Filon-type method requires the availability of the first few moments $\int_{-1}^1x^me^{{\rm i}\omega g(x)}dx$ or the modified moments  $\int_{-1}^1T_me^{{\rm i}\omega g(x)}dx$ in an explicit form. Particularly the reliable and stable computation of the modified moments in the presence of stationary points is for a large part still an open subject and needs further research. For moment-free integration we can refer to the methods presented in \cite{Levin1996,Olver2007,Olver2010}.

\section{Conclusion} 
We restricted our discussion to the important problem
\[\int_{-1}^1f(x)e^{{\rm i}\omega x}dx\] and discussed some advanced but practical numerical methods sharing the same important property: their performance drastically improves as the frequency grows. Although EF and Filon-type methods have different points of departure, their basic underlying principle is the same: whereas a classical Gauss quadrature rule interpolates the whole integrand $f(x)e^{{\rm i}\omega x}$ by a polynomial, they interpolate the function $f$ by a polynomial. Different interpolation nodes can be used and by allowing them to depend on the frequency $\omega$ good asymptotic behaviour is produced. We identified which interpolation nodes for $f$ lead to the best asymptotic order and how this can be combined with good behaviour for smaller frequencies. Interpolating at these nodes along with Chebyshev points ensures convergence. Moreover, since the Chebyshev points are nested, an inexpensive adaptive procedure can be based on comparing the method with $2n$ Chebyshev points with the one with $n$ such points.

\end{document}